# A sequential linear programming (SLP) approach for uncertainty analysis-based data-driven computational mechanics

Mengcheng Huang[1], Chang Liu[1,2], Zongliang Du[1,2*], Shan Tang[1,2], Xu Guo[1,2†]

[1]*State Key Laboratory of Structural Analysis for Industrial Equipment,*
*Department of Engineering Mechanics,*
*International Research Center for Computational Mechanics,*
*Dalian University of Technology, Dalian, 116023, P.R. China*
[2]*Ningbo Institute of Dalian University of Technology, Ningbo, 315016, P.R. China*

**Abstract**

In this article, an efficient sequential linear programming algorithm (SLP) for uncertainty analysis-based data-driven computational mechanics (UA-DDCM) is presented. By assuming that the uncertain constitutive relationship embedded behind the prescribed data set can be characterized through a convex combination of the local data points, the upper and lower bounds of structural responses pertaining to the given data set, which are more valuable for making decisions in engineering design, can be found by solving a sequential of linear programming problems very efficiently. Numerical examples demonstrate the effectiveness of the proposed approach on sparse data set and its robustness with respect to the existence of noise and outliers in the data set.

Corresponding author: *zldu@dlut.edu.cn (Z. Du),
Corresponding author: †guoxu@dlut.edu.cn (X. Guo)

1. **Introduction**

Classical computational mechanics theories and algorithms are always developed based on specific constitutive models. In particular, explicit constitutive relations between state variables (e.g., stress and strain) must be established in advance by interpolating a certain amount of experimental/observational data. Although tremendous achievements have been made by this model-based computation mechanics paradigm, it still suffers from some problems such as inevitable modeling errors and uncertainty, artificial assumptions of the constitutive function/functional forms as well as empirical selections of internal variables etc.

In order to bypass the empirical material modeling step in convention computation mechanics paradigm and eliminate the material modeling empiricism, Kirchdoerfer and Ortiz first proposed the paradigm of data-driven computational mechanics [1]. In this seminal contribution, conservation laws and kinematic relationships are formulated as hard constraints in an optimization problem while material data is used directly to characterize the material behavior instead of constructing explicit constitutive models as in the classical model-based computational mechanics paradigm. Since then, data-driven computational mechanics has received ever-increasing research attention and became an active research direction in the field of computational mechanics. To alleviate the influence of data noise, a clustering analysis based approach has been established by the same authors to enhance the robustness of the DDCM approach against outliers [2]. Alternatively, He and Chen suggested using the information of *k*-nearest neighbors to construct a set of local models for robustly approximating the constitutive manifold with outliers [3]. Kanno also proposed a simple heuristic strategy for data-driven static analysis of truss structures with data involving noise and outliers [4]. In addition, it is revealed that the data-driven paradigm proposed in [1] can be reformulated as a mixed-integer quadratic programming problem whose global optimal solution can be obtained by the branch-and-bound method in principle [5]. Moreover, under some regularity assumptions, Conti et al. proved the existence of solution and the convergence of the corresponding numerical solution approach under the DDCM framework for elasticity problems [6, 7]. In recent years, DDCM approaches have also been generalized from the linear statics analysis to dynamic structural analysis [8], geometric nonlinear analysis [7], [9], [10], diffusion problems [11], fracture modeling [12], anisotropy elasticity [13] and simulation of history-dependent mechanical behaviors [14]. Inspired by the DDCM approach, Leygue et al. proposed a data-driven identification (DDI) algorithm [15, 16], which is capable of identifying the stress field from the measured strain field and the prescribed external loads without resorting to any constitutive model. Later on, this inverse method has also been generalized to allow for elasto-plasticity [17], nonlinear elasticity [18], and elastodynamics problems [19], respectively. Different from the treatment in [1], Ibanez et al. [20, 21] suggested using manifold learning techniques to characterize the material constitutive behavior locally so as to improve the efficiency of data-driven solver. Based on this idea, some effective data-driven approaches have been developed in [22, 23]. In order to deal with high-dimensional (i.e., two dimensional (2D)/ three dimensional (3D)) problems under the DDCM paradigm effectively, a tensor voting approach has been proposed in [24] based on eager machine learning techniques. He et al. developed a deep autocoding technique to learn low-dimensional representations of high-dimensional data sets to improve computing efficiency of the data-driven solution process [25]. For further speeding up the DDCM approach, Eggersmann et al. developed an approximate nearest-neighbor algorithm which can deal with one billion material data points efficiently in high dimensional phase space [26]. Furthermore, considering the fact that it is difficult to obtain experimental data points in high-dimensional phase space required by data-driven approaches, a so-called MAP123 approach has been proposed in [27-29], which can realize efficient data-driven computation using only appropriately selected one-dimensional experiment data.

Although the classical DDCM framework sets up a new model-free paradigm to solve computational mechanics problems by utilizing material data directly, it also faces some challenging issues deserving further explorations. For example, compared with the traditional model-based paradigm, numerical solutions under the DDCM paradigm usually involve more computational effort. This is due to the fact that the mathematical formulation of DDCM is actually a bi-level program, and the lower level program aiming at finding the closest material data corresponding to a specific mechanical state is combinatorial in nature. Since the numbers of the sampling material points/finite elements should be large enough to realize a reliable material property characterization /spatial discretization, and the computational complexity of the searching process for material state is directly proportional to these two numbers, the DDCM-based solution process may be very time consuming especially for 3D problems. Moreover, the bi-level formulation lacks the necessary "differentiable" structure rendering the application of the tools of differential calculus and calculus of variations. This unpleasant behavior may further deteriorate the computation efficiency. Another issue worthy of pointing out is that, as discussed in [30], taking the unavoidable multi-source uncertainties possibly arising from measurement errors, information deficiency and model inaccuracy associated with the material data collection process into consideration, it seems reasonable to consider *a solution set* rather than *a single solution* when the data-drive paradigm is employed for solving computational mechanics problems. This is because the former one with confidence bounds is more helpful for making decisions in engineering applications than the latter one when the existence of uncertainties in material property characterization is inevitable.

Based on the aforementioned considerations, a new uncertainty analysis-based data-driven computational mechanics (UA-DDCM) framework was proposed in [30]. Compared with the original DDCM framework, the UA-DDCM framework focuses on obtaining *a confidence bounding interval* of a concerned structural response rather than *a single nominal value*. To this end, it was proposed to cover the data set tightly by a set of ellipsoids and formulate the corresponding UA-DDCM problem as *a single level* mathematical program. It has also been shown that when the data set can be enclosed tightly by a single ellipsoid, the corresponding problem will be *convex* in nature and can be solved with very efficient algorithms. It is worth noting, however, that even though the problem can be transformed into a single level program under the UA-DDCM framework, the corresponding problem is in general non-convex and non-smooth when multiple ellipsoids are employed to cover the data set. The corresponding solution procedures are also not discussed in detail in [30]. In the present work, in order to solve the single level program formulated under the UA-DDCM framework efficiently, a sequential linear programming (SLP) approach is developed. The central idea is to solve a sequence of linear programming problems by constructing a local convex hull of a number of data points in the constitutive data set for each structural member (for discrete case)/ Gauss point (for continuum case), adaptively. It is assumed that the points locating in the convex hull (even though not coinciding with any data point) may also represent the constitutive behaviors of the considered material. This treatment in some sense characterizes the possible uncertainty embedded behind the prescribed data set usually obtained from a limited number of physical experiments and/or numerical simulations. The sizes of the convex hulls are first set to take some relatively large values and then reduced gradually during the process of sequential optimization following a trust-region like strategy. It is found through numerical experiments that the proposed approach can not only enhance the performance of the original DDCM-based algorithms when the prescribed data set is only comprised of a limited number of data points and contains noise and outliers, but also provide a quantitative measure of the influence of uncertainties on concerned structural responses.

The rest of this article is organized as follows. In Section 2, the formulation of the classical data-driven computational mechanics (DDCM) framework and the so-called uncertainty analysis-based data-driven computational mechanics (UA-DDCM) framework, which can account for the uncertainty of constitutive relationship embedded behind the prescribed data set, are described briefly. Then a sequential linear programming (SLP) approach developed under the UA-DDCM framework is introduced in detail in Section 3. In Section 4, the key idea of the proposed algorithm is first illustrated by a two-dimensional truss example and then the convergence property, robustness, accuracy and the importance of considering uncertainties in data-driven computational mechanics are verified by solving a three-dimensional truss structure. Finally, a cantilever beam example is examined to illustrate the potential of the SLP-UADDCM algorithm for tackling three-dimensional continuum problems. At the end of this article, some concluding remarks and perspectives for future research works are presented.

## 2. The data-driven computational mechanics (DDCM) frameworks

In this section, for the sake of completeness, the classical DDCM framework proposed by Kirchdoerfer and Ortiz [1] and the UA-DDCM framework capable of accounting for the uncertainty of data based characterization of constitutive relationship proposed in [30] are briefly described.

**2.1 The classical data-driven computational mechanics (DDCM) framework**

The general formulation for analyzing elastic structure in the classical DDCM framework can be found in [1]. For the ease of illustration, the corresponding mathematical formulation for truss structures in the classical DDCM framework is reviewed at first. As demonstrated in [5], under the assumption of infinitesimal deformation and linear elasticity, structural analysis of a truss structure in the DDCM framework can be formulated as the following optimization problem:

$$\text{Find} \quad \boldsymbol{d} = (\boldsymbol{U}^\top, \boldsymbol{\sigma}^\top, \boldsymbol{\varepsilon}^\top, \boldsymbol{t}^\top)^\top$$

$$\text{Min} \quad I = \sum_{e=1}^{m} \frac{1}{2} v_e c_e (\varepsilon_e - \epsilon_e)^2 + \sum_{e=1}^{m} \frac{1}{2} \frac{v_e}{c_e} (\sigma_e - s_e)^2$$

$$\text{S.t.} \quad \varepsilon_e = \boldsymbol{b}_e^\top \boldsymbol{U} / l_e, \quad e = 1, \dots, m,$$

$$\sum_{e=1}^{m} A_e \sigma_e \boldsymbol{b}_e = \boldsymbol{p}, \quad e = 1, \dots, m,$$

$$\begin{pmatrix} \epsilon_e \\ s_e \end{pmatrix} = \sum_{j=1}^{N_d} \begin{pmatrix} \varepsilon_j^d \\ \sigma_j^d \end{pmatrix} t_{ej}, \quad e = 1, \dots, m,$$

$$\sum_{j=1}^{N_d} t_{ej} = 1, \quad e = 1, \dots, m,$$

$$t_{ej} \in \{0, 1\}, \quad e = 1, \dots, m; \ j = 1, \dots, N_d, \tag{1}$$

where $\boldsymbol{U} = (U_1, \dots, U_n)^\top$ and $\boldsymbol{p} = (p_1, \dots, p_n)^\top \in \mathbf{R}^n$ are the vectors of nodal displacement and external load with $n$ denoting the number of degree-of-freedom, the symbols $\boldsymbol{\sigma} = (\sigma_1, \dots, \sigma_m)^\top$ and $\boldsymbol{\varepsilon} = (\varepsilon_1, \dots, \varepsilon_m)^\top \in \mathbf{R}^m$ are the vectors of the stresses and strains of truss bars with $m$ denoting the total number of bars in the truss structure. The vector $\boldsymbol{t} = (t_{11}, \dots, t_{eN_d}, \dots, t_{m1}, \dots, t_{mN_d})^\top \in \mathbf{R}^{m \times N_d}$ with only binary components is used for identifying a specific data point in the given data set $\mathcal{D} = \{(\varepsilon_1^d, \sigma_1^d), \dots, (\varepsilon_{N_d}^d, \sigma_{N_d}^d)\}$ ($N_d$ denotes the total number of data points). Besides, the quantities $v_e$, $A_e$, and $l_e$ are the volume, the cross-sectional area and the length of the $e$-th bar, respectively. In addition, $c_e$ is a scaling factor and the symbol $\boldsymbol{b}_e$ is the vector of the director cosine of the $e$-th bar.

In the classical DDCM framework, it is intended to find *a single solution* which, besides satisfying the conservation laws and compatibility conditions, has the closest distance to a prescribed data set $\mathcal{D}$ characterizing the material behavior in the phase space. The most distinctive feature in the DDCM is that there is no need to establish an explicit constitutive model and the constitutive relationship is preserved point-wisely (or element-wisely) through a data-driven distance-minimizing scheme [1]. Although the classical DDCM framework opens a new avenue for computational mechanics, as an emerging field, some challenging issues described in the introduction have been undergoing intensive explorations since its invention.

**2.2 The uncertainty analysis-based data-driven computational mechanics (UA-DDCM) framework**

In order to account for the influence of the unavoidable multi-source uncertainties in the data set on data-driven solutions, Guo et al. [30] developed a uncertainty analysis-based framework for data-driven computational mechanics. In this framework, it is proposed to enclose the prescribed data set with possible outliers by a union of totally $L$ ellipsoids in the stress-strain space with a minimum volume (as illustrated in Fig. 1c). Furthermore, it is also assumed that all the points (not only the prescribed experimental data points!) inside the enclosed ellipsoids may represent the possible constitutive behavior of the considered material. Therefore, instead of pursuing a single solution as in the classical DDCM approach, *a solution set*, which includes the extreme values of the concerned structure response, should be determined. Under this consideration, the UA-DDCM framework for structural analysis of truss structures can be formulated as follows:

$$\begin{aligned}
&\text{Find} \quad \boldsymbol{d} = (\boldsymbol{U}^\top, \boldsymbol{\sigma}^\top, \boldsymbol{\varepsilon}^\top)^\top \\
&\text{Min} \quad I(\boldsymbol{d}) = \boldsymbol{s}^\top \boldsymbol{d} \\
&\text{S.t.} \quad \varepsilon_e = \boldsymbol{b}_e^\top \boldsymbol{U}/l_e, \quad e = 1, \dots, m, \\
&\qquad \sum_{e=1}^{m} A_e\, \sigma_e \boldsymbol{b}_e = \boldsymbol{p}, \quad e = 1, \dots, m, \\
&\qquad \min\bigl(h_1(\varepsilon, \sigma), \dots, h_L(\varepsilon, \sigma)\bigr) \leq 0,
\end{aligned} \qquad (2)$$

where $\boldsymbol{s} = (0,0,\dots,\pm 1_{i\text{th}},0,\dots,0)^\top \in \mathbf{R}^{n+2m}$ is an indication vector. In Eq. (2), $h_j,\ j = 1,\dots,L$ has following form:

$$h_j(\varepsilon, \sigma) = \left(\varepsilon - \varepsilon_j^0, \sigma - \sigma_j^0\right) \boldsymbol{P}_j \left(\varepsilon - \varepsilon_j^0, \sigma - \sigma_j^0\right)^\top - 1, \qquad (3)$$

where $\left(\varepsilon_j^0, \sigma_j^0\right)$ and $\boldsymbol{P}_j$ denote the center point and the shape matrix of the $j$-th ellipsoid, respectively.

It is worth noting that taking the uncertainty embedded in the data set into consideration, the above UA-DDCM framework has the potential of providing a confidence bound of the concerned structural response. This is very important for making decisions in practical engineering applications compared to the case when only a single nominal value is available. The UA-DDCM framework is also robust with respect to the existence of outliers in the data set since the corresponding mathematical formulation renders the possibility of searching the data set in a global way [31]. Moreover, since Eq. (3) is actually a single-level program with continuous variables, the corresponding solution process is theoretically more efficient than that of the classical DDCM framework which is a mixed 0-1 program involving both continuous and discrete variables in nature.

Although the UA-DDCM framework was established in [30], the corresponding solution procedure was not discussed in detail. In the following section, a sequential linear programming-based approach is proposed to address this issue. The central idea is to utilize the convexity property embedded in the underlying problem to enhance the efficiency and robustness of the solution process.

## 3. A sequential linear programming approach under the UA-DDCM framework

### 3.1 Problem formulation

Instead of constructing a set of ellipsoids to cover the data set, it is also possible to bound the data set by a single polygon as shown in Fig. 1d. This treatment can not only provide a tighter encloser of the data set, but also transfer the non-convex problem in Eq. (2) including a set of quadratic constraints into a convex one with a number of *linear* constraints, which can be solved by very powerful modern linear programming approaches (e. g., the interior point type algorithms).

Although constructing a single polygon can enhance the efficiency of finding the bounds of structural responses significantly under the UA-DDCM framework, the data set will be over-relaxed under this treatment especially when the constitutive behavior represented by the data points is far from the "linear" form or numerous outliers exist in the data set. Under these circumstances, the gap between the obtained upper and lower bounds may be very large. This inspires us to construct the convex hulls locally and update them iteratively to approach the "true" constitutive responses of the involved material represented by the appropriate sets of the prescribed data points. Based on this consideration, a sequential linear programming formulation can be constructed as follows (also taking the truss structure as an example) $\mathcal{P}^{(k)}$:

$$\text{Find} \quad \boldsymbol{d}^{(k)} = \left((\boldsymbol{U}^{(k)})^\top, (\boldsymbol{\varepsilon}^{(k)})^\top, (\boldsymbol{\sigma}^{(k)})^\top\right)^\top$$

$$\text{Min} \quad I(\boldsymbol{d}) = \boldsymbol{s}^\top \boldsymbol{d}^{(k)}$$

$$\text{S.t.} \quad \varepsilon_e^{(k)} = \boldsymbol{b}_e^\top \boldsymbol{U}^{(k)}/l_e, \quad e = 1, \dots, m,$$

$$\sum_{e=1}^m A_e \sigma_e^{(k)} \boldsymbol{b}_e = \boldsymbol{p}, \quad e = 1, \dots, m,$$

$$\begin{pmatrix} \varepsilon_e^{(k)} \\ \sigma_e^{(k)} \end{pmatrix} = \sum_{j=1}^{N_c} \lambda_{ej}^{(k)} \begin{pmatrix} \left(\varepsilon_{j,e}^{\mathrm{d}}\right)^{(k)} \\ \left(\sigma_{j,e}^{\mathrm{d}}\right)^{(k)} \end{pmatrix}, \quad e = 1, \dots, m,$$

$$\sum_{j=1}^{N_c} \lambda_{ej}^{(k)} = 1, \quad e = 1, \dots, m,$$

$$0 \le \lambda_{ej}^{(k)} \le 1, \quad e = 1, \dots, m;\ j = 1, \dots, N_c. \quad (4)$$

In Eq. (4), $k$ denotes the number of iteration and $\boldsymbol{d}^{(k)}$ represents the value of $\boldsymbol{d}$ to be found in the $k$-th iteration. The symbol $\boldsymbol{\lambda}_e^{(k)} = \left(\lambda_{e1}^{(k)}, \dots, \lambda_{eN_c}^{(k)}\right)^\top$, $e = 1, \dots, m$ is the vector of the coefficients of the convex combination of $N_c$ data points $\left\{\left((\varepsilon_e^{\mathrm{d}})^{(k)}, (\sigma_e^{\mathrm{d}})^{(k)}\right)^\top\right\} = \left\{\left((\varepsilon_{1,e}^{\mathrm{d}})^{(k)}, (\sigma_{1,e}^{\mathrm{d}})^{(k)}\right)^\top, \dots, \left((\varepsilon_{N_c,e}^{\mathrm{d}})^{(k)}, (\sigma_{N_c,e}^{\mathrm{d}})^{(k)}\right)^\top\right\}$ associated with the $e$-th bar of the truss structure. Here the involved sets of data points $\left\{\left((\varepsilon_e^{\mathrm{d}})^{(k)}, (\sigma_e^{\mathrm{d}})^{(k)}\right)^\top\right\}$, $e = 1, \dots, m$ must be specified in advance when $\mathcal{P}^{(k)}$ is solved. These data sets, however, will be updated adaptively during the course of sequential iteration in a way described in the subsequent text.

## 3.2 Solution procedure

In the proposed approach, the upper/lower bound of a concerned structural response pertaining to a given constitutive data set is found by solving a series of linear programs ($\mathcal{P}^{(k)}$ in Eq. (4)) in a sequential way. In the following, the details of the corresponding numerical implementation will be described. The outline of the proposed sequential linear programming (SLP) approach for the UA-DDCM framework is summarized in Table 1. To be specific, the proposed algorithm can be decomposed into the following four steps:

(1) Initialization of the local data set for convex hull construction

To guarantee the feasibility of the linear programming in (4), the initial convex hull for identifying the constitutive behavior should be sufficiently large. In our implementation, all the data points in $\mathcal{D}$ are first sorted according to the values of $\text{sign}(\varepsilon^\text{d})\|(\varepsilon^\text{d}, \sigma^\text{d})\|_2$. Then without loss of generality, the initial data points first used for the convex hull construction of each bar can be uniformly chosen as: $\left\{(\varepsilon_1^\text{d}, \sigma_1^\text{d}), (\varepsilon_{1+L^{(1)}}^\text{d}, \sigma_{1+L^{(1)}}^\text{d}), \dots, (\varepsilon_{1+(N_\text{c}-1)L^{(1)}}^\text{d}, \sigma_{1+(N_\text{c}-1)L^{(1)}}^\text{d})\right\}$ where $L^{(1)}$ is the integer part of $(N_\text{d} - 1)/N_\text{c}$.

(2) Determination of the structural responses

If the linear programming $\mathcal{P}^{(k)}$ in Eq. (4) is feasible (this is mostly often the case from our numerical experience), structural responses such as $\boldsymbol{U}^{(k)}, \boldsymbol{\varepsilon}^{(k)}$ and $\boldsymbol{\sigma}^{(k)}$ can be calculated directly by very efficient algorithms (e.g., simplex or interior point algorithm). If, however, $\mathcal{P}^{(k)}$ is not feasible (this is actually rarely encountered in our numerical experiments), we first calculate the geometric center of $\left((\varepsilon_{j,e}^\text{d})^{(k)}, (\sigma_{j,e}^\text{d})^{(k)}\right)$, $j = 1, \dots, N_\text{c}$; $e = 1, \dots, m$ as:

$$\left(\tilde{\varepsilon}_e^{(k)}, \tilde{\sigma}_e^{(k)}\right) = \frac{1}{N_\text{c}} \left( \sum_{j=1}^{N_\text{c}} (\varepsilon_{j,e}^\text{d})^{(k)}, \sum_{j=1}^{N_\text{c}} (\sigma_{j,e}^\text{d})^{(k)} \right), \qquad e = 1, \dots, m. \qquad (5)$$

Once $\left(\tilde{\varepsilon}_e^{(k)}, \tilde{\sigma}_e^{(k)}\right)$, $e = 1, \dots, m$ is determined, the corresponding structural response $\boldsymbol{U}^{(k)}$ and an intermediate multiplier $\boldsymbol{\eta}^{(k)}$ can be calculated by solving the following equations firstly [1]:

$$\sum_{j=1}^{n}\left(\sum_{e=1}^{m}A_e b_{ej} b_{ei}/l_e^2\right)U_j^{(k)} = \sum_{e=1}^{m}A_e b_{ei}\tilde{\varepsilon}_e^{(k)}/l_e, \tag{6a}$$

$$\sum_{j=1}^{n}\left(\sum_{e=1}^{m}A_e b_{ej} b_{ei}\right)\eta_j^{(k)} = p_i - \sum_{e=1}^{m}A_e b_{ei}\tilde{\sigma}_e^{(k)} \tag{6b}$$

and then obtain $\boldsymbol{\varepsilon}^{(k)}$ and $\boldsymbol{\sigma}^{(k)}$ as:

$$\varepsilon_e^{(k)} = \sum_{i=1}^{n} b_{ei} U_i^{(k)}/l_e, \quad \sigma_e^{(k)} = \tilde{\sigma}_e^{(k)} + \sum_{i=1}^{n} b_{ei}\eta_i^{(k)}, \quad e = 1, \dots, m, \tag{7}$$

respectively.

(3) Update of the data points for local convex hull construction

Based on the strategies for obtaining nodal displacement, strain and stress vectors, the following two cases are considered for updating the local data sets for the convex hull construction.

Case 1: The linear programming $\mathcal{P}^{(k)}$ in Eq. (4) is feasible, but the local convex hull in stress-strain space may be over-sized in the $k$-th iteration, as shown in Fig. 1d. Under this circumstance, in order to obtain a tighter bound, the size of the convex hull in the $(k + 1)$-th iteration should be reduced. In our implementation, the Euclidean distances between all data points and $\left(\varepsilon_e^{(k)}, \sigma_e^{(k)}\right)$, $e = 1, \dots, m$, i.e., $d\left(\left(\varepsilon_e^{(k)}, \sigma_e^{(k)}\right), \left(\varepsilon_j^{\mathrm{d}}, \sigma_j^{\mathrm{d}}\right)\right)$, $j = 1, \dots, N_{\mathrm{d}}$ will be calculated first. Then the index $ID_e^{(k)}$ for the data point which is closest to $\left(\varepsilon_e^{(k)}, \sigma_e^{(k)}\right)$ for every $e = 1, \dots, m$ can be determined in $\mathcal{D}$ based on the values of $d\left(\left(\varepsilon_e^{(k)}, \sigma_e^{(k)}\right), \left(\varepsilon_j^{\mathrm{d}}, \sigma_j^{\mathrm{d}}\right)\right)$, $j = 1, \dots, N_{\mathrm{d}}$. Since the data points in $\mathcal{D}$ have already been sorted according to the values of $\mathrm{sign}(\varepsilon^{\mathrm{d}})\|(\varepsilon^{\mathrm{d}}, \sigma^{\mathrm{d}})\|_2$ when the local data set is initialized as discussed in step (1), we only need to introduce an integer indicator $L^{(k)} \geq 1$ to measure the size of the convex hull in the $k$-th iteration, and update $L^{(k+1)}$ as the integer part of $L^{(k)}/\rho$ with $\rho > 1$. Therefore, the local data points involved in the convex hull construction for each bar ($e = 1, \dots, m$) in the $(k + 1)$-th iteration can be selected as:

$$\left\{\left(\varepsilon_{ID_e^{(k)}-tL^{(k+1)}}^{\mathrm{d}}, \sigma_{ID_e^{(k)}-tL^{(k+1)}}^{\mathrm{d}}\right), \dots, \left(\varepsilon_{ID_e^{(k)}}^{\mathrm{d}}, \sigma_{ID_e^{(k)}}^{\mathrm{d}}\right), \dots, \left(\varepsilon_{ID_e^{(k)}+tL^{(k+1)}}^{\mathrm{d}}, \sigma_{ID_e^{(k)}+tL^{(k+1)}}^{\mathrm{d}}\right)\right\},$$

with $t$ denoting the integer part of $N_{\mathrm{c}}/2$.

Case 2: The linear program $\mathcal{P}^{(k)}$ in Eq. (4) is infeasible. In this case, the current local convex hull needs to be enlarged to guarantee the feasibility of $\mathcal{P}^{(k)}$. Therefore, instead of reducing the size of the local convex hull, we suggest to use a larger $L^{(k+1)}$ ($L^{(k+1)} = L^{(k)} + 1$ in the present work) to determine the local data points for the $(k+1)$-th iteration following the procedure described in Case 1.

(4) Check of convergence

The iteration process terminates once the relative error of the $L_2$-norm of the displacement vector is less than a threshold value.

## 4. Numerical examples

In this section, a set of examples are investigated and discussed to evaluate the performance of the proposed approach. A two dimensional three-bar truss example is first used to illustrate the key ideas of the present SLP-UADDCM algorithm. Afterwards, a three-dimensional truss structure and a cantilever beam structure discretized by finite elements are analyzed to examine the robustness, effectiveness and accuracy of the proposed approach. All examples are solved on a laptop equipped with an Intel(R) Core(TM) 2.61GHz CPU and 32.0GB of RAM.

**4.1 A three-bar truss example**

In this example, shown in Fig. 2a, a 2D truss structure with $A = 1$ for all bars, $l_1 = l_3 = 1$, $l_2 = \sqrt{2}$ and $|\boldsymbol{p}| = \sqrt{2}/2$ is considered. As shown in Fig. 2b, the deliberately designed noisy data set including 201 points is generated based on a reference linear elastic constitutive relationship as $\sigma_j^{\mathrm{d}} = E\varepsilon_j^{\mathrm{d}} - \vartheta + 2\vartheta\mathcal{U}(0,1)$, $j = 1, \ldots, 201$ with $E = 1$ and $\vartheta = \left|E\varepsilon_j^{\mathrm{d}}\right|$ for $\left|E\varepsilon_j^{\mathrm{d}}\right| \leq 0.1$ ($\vartheta = 0.1$ otherwise). Here, $\mathcal{U}(0,1)$ is a random value distributed uniformly in $[0,1]$. The solution procedure described in Table 1 is adopted to solve this problem by setting $L^{(1)} = 25$, $\rho = 1.5$, $N_c = 5$ and $\mathrm{Tol} = 0.01$, respectively.

Since the exact stress and strain states in each bar are generally unknown initially, in order to guarantee the feasibility of the linear programming $\mathcal{P}^{(1)}$, a sufficiently large convex hull should be constructed for each bar in the first iteration. This is achieved by selecting data points with relatively large distances in the data set for the convex hull construction (see the five data points colored in red in Fig. 2c). Of course, closer local data points can be chosen if the initial strain/stress states of some bars can be estimated a priori to accelerate the convergence of the iteration process. Numerical experiments indicate that although the local data points chosen for convex hull construction are the same for the three bars initially, different local convex hulls can be identified efficiently and updated adaptively based on the values of the stresses and strains of the bars obtained in the previous iteration step as shown in Fig. 2c.

To obtain the upper and lower bounds of the horizontal displacement of the free node, the objective function can be set as $I(\boldsymbol{U}) = -U_1$ and $I(\boldsymbol{U}) = U_1$, respectively. The variations of the values of the horizontal displacement ($U_1$) and the vertical displacement ($U_2$) of the free node during the process of iteration obtained with different objective functions are listed in Table 2. It is observed from Fig. 2c that as the local convex hulls are gradually shrunk, the upper bound of $U_1$ (i.e., $\overline{U}_1$ obtained by setting $I(\boldsymbol{U}) = -U_1$) decreases from 0.5856 in the first iteration to 0.5334 at the 12-th iteration. Similarly, for the case $I(\boldsymbol{U}) = U_1$, the lower bound of $U_1$ (i.e., $\underline{U}_1$) increases from 0.3764 in the first iteration to 0.4184 at iteration 10 (see Fig. 2d for the evolution of the local data points for convex hull construction). Accordingly, the gap between $\overline{U}_1$ and $\underline{U}_1$ drops from 0.2091 to 0.1150, which implies that the confidence of uncertainty quantification improves significantly through the optimization process. In addition, the linear program in (5) was also solved by setting the objective function as $\boldsymbol{p}^\top \boldsymbol{U}$ and the corresponding converged local data points are plotted in Fig. 2e. As shown in Table 2, when $I(\boldsymbol{U}) = \boldsymbol{p}^\top \boldsymbol{U}$ is adopted, the converged value is $U_1 = 0.4226$, which consistently falls into the interval of [0.4184, 0.5334] determined by solving the aforementioned optimization problems with $I(\boldsymbol{U}) = \pm U_1$ separately. These results clearly demonstrate the effectiveness of the UA-DDCM formulation and the present approach for its numerical implementation.

In addition, the noisy data can be approximately enveloped by two linear elastic constitutive relations with $E = 0.8$ and $E = 1.2$, respectively. From this point of view, this problem can also be solved as an extreme analysis of truss structures with material uncertainties [32]. Correspondingly, the response interval of $U_1$ is determined as [0.4167, 0.625], which is much larger than the response bounds of $[0.4184, 0.5334]$ determined by the proposed SLP-UADDCM algorithm. This is because in the proposed algorithm, the uncertainty sets, i.e., the local convex hulls, are updated adaptively during the iteration process and much smaller than the fixed uncertainty set in classical robust optimization method. Furthermore, in order to validate the ability of the proposed algorithm for obtaining theoretical bounds, a regularized data set without noise shown in Fig. 3a is examined. The convergence histories of $U_1$ obtained by the SLP-UADDCM algorithm with $I(U) = \pm U_1$ are shown in Fig. 3b and the obtained bounds are exactly the same as the theoretical values (i.e., $\underline{U}_1 = 0.4167$ and $\overline{U}_1 = 0.6250$, respectively).

**4.2 A 3D truss example**

In this example, the convergence and robustness properties of the SLP-UADDCM algorithm are verified by analyzing a three-dimensional truss structure with 1194 bars (1002 degrees of freedom) described in Fig. 4a. It is assumed the truss undergoes small deformation and the material composed of the bars obeys a nonlinear elastic constitutive law (i.e., $\sigma = \sigma(\varepsilon) = \varepsilon^{1/3}$) shown in Fig. 4b. A Newton-Raphson solver is employed to obtain the reference model-based solutions for comparisons.

*4.2.1 Solutions with objective function $I(U) = p^\mathrm{T} U$*

In this subsection, the objective function is chosen as $I(U) = p^\mathrm{T} U$ and the corresponding solutions are compared with the counterparts obtained under the classical computational mechanics framework.

1) Convergence of the SLP-UADDCM algorithm

The convergence property of the proposed algorithm is evaluated with a precise data set from the following aspects: the evolution of the local data points for convex hull construction, the variation of the concerned structural response during the iteration process and the influence of the size of the data set used in the algorithm on the convergence process. Related parameters in Table 1 are set as $L^{(1)} = 25$, $\rho = 2$, $N_c = 5$ and $\mathrm{Tol} = 0.001$, respectively.

By generating 121 data points following the exact constitutive relationship $\sigma = \varepsilon^{1/3}$, a converged solution is obtained by the proposed approach within 6 iterations. Fig. 5 provides the evolution of the local data points used for convex hull construction associated with the 884-th bar element (see Fig. 4a for reference). It is observed that the size of the local data set is gradually reduced as the iteration process proceeds. Furthermore, the relative error of displacement vector and the normalized root-mean-square (RMS) errors of stress and strain vectors are evaluated as:

$$U_{\text{RE}} = \frac{\|U^{(k)} - U^{\text{ref}}\|_2}{\|U^{\text{ref}}\|_2}, \quad \sigma_{\text{RMS}} = \frac{\|\sigma^{(k)} - \sigma^{\text{ref}}\|_2}{\sqrt{m}\|\sigma^{\text{ref}}\|_\infty} \quad \text{and} \quad \varepsilon_{\text{RMS}} = \frac{\|\varepsilon^{(k)} - \varepsilon^{\text{ref}}\|_2}{\sqrt{m}\|\varepsilon^{\text{ref}}\|_\infty},$$

respectively, where $U^{\text{ref}}$, $\sigma^{\text{ref}}$ and $\varepsilon^{\text{ref}}$ are the reference displacement, stress and strain vectors obtained by the conventional Newton-Raphson algorithm with $\sigma = \varepsilon^{1/3}$. The corresponding iteration process of $U_{\text{RE}}$ and $\sigma_{\text{RMS}}$ is plotted in Fig. 6, which reveals that the structural response obtained using $I(U) = p^\top U$ as the objective function is quite close to the reference solution in this example when the exact constitutive data set is employed.

Next, the convergence behavior of the proposed algorithm with respect to the total number of data points is investigated. For all data sets examined, we choose $\rho = 2$, $N_c = 5$, and $L^{(1)}$ as the integer part of $N_d/N_c$ plus 1. The values of $U_{\text{RE}}$, $\sigma_{\text{RMS}}$ and $\varepsilon_{\text{RMS}}$, the number of iterations for convergence, and time costs for different tested cases are presented in Table 3. The computational results clearly demonstrate that as the data set approaches the exact constitutive model, all of $U_{\text{RE}}$, $\sigma_{\text{RMS}}$ and $\varepsilon_{\text{RMS}}$ decrease and the corresponding data-driven result converges to the reference model-based solution. Another interesting point is that the converged iteration numbers and the solution times of the proposed algorithm do not increase significantly as the total number of data points increasing from 41 to 100001. This can be contributed to the updating strategy of local data points described in Table 1. Actually, numerical experiments indicate that a larger value of $\rho$ could further increase the convergence rate.

It is also worth noting that even for this relatively complex 3D truss structure composed by nonlinear material, the solution time is only about 1s-2s on a laptop, which is close to the time cost of the classical Newton-Raphson solver (1.43s). This implies that the proposed data-driven algorithm may also find its application even under the model-based solution framework by representing the explicit constitutive function using a set of discrete data points with very close distances.

In addition, as illustrated in Fig. 7, the solutions obtained by the proposed approach converge to the reference solution as the total number of data points increases. Compared with the convergence results presented in [1], the present approach has smaller RMS errors of stress and strain when the number of data points is less than $10^5$. This can be understood from that, since all the stress-strain pairs inside the local convex hull are feasible, on the one hand, when the data points are insufficient, this treatment could enrich the data set effectively; on the other hand, for the case there are sufficiently dense data points, the local convex hull still may introduce stress-strain pairs not exactly locating on the constitutive manifold even though there is no noise on the data points for convex hull construction. This character would be attractive when the data points are not easy to be obtained or the curse of dimensionality for three dimensional problems exists.

2) Robustness of the SLP-UADDCM algorithm

Artificial noises are deliberately added to the constitutive curve plotted in Fig. 4b according to the relation $\sigma_j^د = \left(\varepsilon_j^د\right)^{1/3} - \vartheta + 2\vartheta\mathcal{U}(0,1)$, $j = 1, \ldots, N_د$, with $N_د = 121$, and $\vartheta = \vartheta_0$ when $\vartheta_0 \leq \left|\left(\varepsilon_j^د\right)^{1/3}\right|$ ($\vartheta = \left|\left(\varepsilon_j^د\right)^{1/3}\right|$ otherwise), respectively. Fig. 8 illustrates a random data set with $\vartheta_0 = 0.04$. In order to fully validate the robustness of the SLP-UADDCM algorithm, we generate three data groups (each contains 100 data sets) with $\vartheta_0 = 0.02, 0.04, 0.08$, respectively. By setting parameters in Table 1 as $L^{(1)} = 25$, $\rho = 1.1$, $N_c = 5$ and $\text{Tol} = 0.01$, the mean and variance of $U_{\text{RE}}$ and $\sigma_{\text{RMS}}$ are presented in Table 4 for different values of $\vartheta_0$, respectively, and Fig. 9 further shows their distribution histograms. It is evident that as the randomness decreases, the corresponding mean values of $U_{\text{RE}}$ and $\sigma_{\text{RMS}}$ get smaller and this demonstrates the robustness of the proposed algorithm.

To further test the performance of the SLP-UADDCM algorithm about outliers, the stress amplitudes of 4, 8, 16, 32 random data points selected from the data set ($\vartheta_0 = 0.04$, $N_d = 121$) are scaled by 0.8 and 1.2 times, respectively. Fig. 10 illustrates the two representative data sets with 32 outliers. With the same parameters setting, the obtained mean and variance values of $U_{\text{RE}}$ and $\sigma_{\text{RMS}}$ are shown in Table 5 and Table 6 for different numbers of outliers and the two scale factors, respectively. The corresponding distribution histograms of $U_{\text{RE}}$ and $\sigma_{\text{RMS}}$ are plotted in Figs. 11a-b, respectively.

According to the variance values listed in Table 5 and Table 6, the algorithm achieves good robustness for different numbers of outliers. Meanwhile, it is also found that the number of outliers has little influence on the obtained stresses while its influence on displacement is more significant. This can be understood from the fact that the stresses and strains are located in the local convex hulls and need to satisfy the equilibrium equation. The displacement vector, however, only needs to satisfy the compatibility conditions and the displacement constraints, and thus has a wider range of values. Similar tendencies can also be found in the histograms in Figs. 10-11, which imply that the noise and outliers have more considerable effect on the values of nodal displacements. Furthermore, the results in Table 5 and Table 6 indicate that outliers scaled by 1.2 times generally have a more significant influence than outliers scaled by 0.8 times on the results. This is because for the adopted objective function $I = \boldsymbol{p}^\top \boldsymbol{U}$, a stronger material would decrease the structural compliance and thus the optimization algorithm tends to select the "stronger material" (i.e., the outliers scaled by 1.2 times) to resist the external load.

*4.2.2 Upper and lower bounds with objective function $I(\boldsymbol{U}) = \pm U_i$*

As shown by the results in subsection 4.2.1, a single solution obtained in the original data-driven framework would be inevitably affected by the noise and outliers in the data sets. In particular, for the case $\vartheta_0 = 0.08$, the relative error of displacement vector is more than 13%, which cannot be neglected in practical engineering applications. Actually, it could be very difficult to identify and exclude outliers or noise from the data set in real engineering applications. Therefore, it would be more practical to present the bounds of the concerned response to allow for the inevitable uncertainties associated with the data set. In this section, the objective function is changed to $\pm U_i$ to obtain the upper and lower bounds of displacement at $i$-th degree of freedom.

1) Convergence of the SLP-UADDCM algorithm

Using the same parameters setting in Table 1, Table 7 presents the upper and lower bounds of $U_{552}$ (the $z$-directional displacement of node 184 in Fig. 4a). Notably, the uncertainty interval of $U_{552}$ drastically decreased from $[-2.1459, -1.4727]$ (with totally 61 data points) to $[-1.6473, -1.6394]$ (with totally 2363 data points), while the reference solution is $U_{552}^{\text{ref}} = -1.6415$. This clearly shows that, for precise data sets, the present SLP-UADDCM approach can not only estimate a practical response bound, but also does have the ability of converging to the exact response as the data set approaching the constitutive manifold.

2) Robustness of the SLP-UADDCM algorithm

In order to explore the influence of noise and outliers, the comparable solution (obtained by setting $I(\boldsymbol{U}) = \boldsymbol{p}^\top \boldsymbol{U}$) as well as the upper and lower bounds (with $I(\boldsymbol{U}) = \pm U_{552}$) of 100 random data sets ($\vartheta_0 = 0.04, N_d = 121$) are shown in Fig. 12a. In those cases, the ranges described by upper bounds and lower bounds always cover both the reference value ($U_{552}^{\text{ref}} = -1.6415$) and the corresponding comparable solutions. Furthermore, since $U_{552}$ is negative, its upper bound corresponds to a smaller deformation and a smaller structural compliance. This explains the phenomenon that the upper bounds are always closer to the corresponding comparable values than the lower bounds. We also present the upper and lower bounds of $U_{70}$ (the $x$-directional displacement of node 24 in Fig. 4a) with $U_{70}^{\text{ref}} = 0.2061$ using the same 100 random noisy data sets shown in Fig. 12b. This time the comparable solution is closer to the corresponding lower bounds as expected.

In order to further explore the robustness of the SLP-UADDCM algorithm for calculating upper and lower bounds, 16 random outliers either scaling the stress amplitudes by 1.2 or 0.8 times are added to the above noisy data sets, as shown in Fig. 10. Fig. 13 shows the corresponding results of $U_{552}$ in the 100 data sets with outliers. Compared with Fig. 12a, by introducing outliers scaled by 1.2 times, both the upper bounds and the comparable solutions in Fig. 13a increase more significantly as compared to the corresponding lower bounds. This is because the outliers correspond to the response of "stronger material", which is effective to decrease the amplitude of deformation and the structural compliance. Based on the same reason, the outliers scaled by 0.8 times decrease the lower bounds more significantly while having little effect on the upper bounds and comparable solutions as illustrated in Fig. 13b.

Based on the above results, it is reasonable to conclude that the proposed SLP-UADDCM framework could present relatively tighter bounds considering the unavoidable uncertainties in the data set and has the capability of covering the reference solution. For the clean data sets without noise, the SLP-UADDCM algorithm can efficiently obtain the exact response as when the number of data points increase. The upper and lower bounds of concerned response are also robust to noise and outliers in data sets, and this not only improves the practical significance but also bypasses some numerical difficulties as compared to the single solution obtained in the classical DDCM framework.

**4.3 A three-dimensional cantilever beam example**

In this section, a dimensionless three-dimensional cantilever beam illustrated in Fig. 14 is studied. This beam is discretized by $16 \times 8 \times 4$ uniform eight-node brick elements, of which the stiffness matrix is calculated using the second-order Gaussian integration. The reference constitutive model used for generation of data sets is a linear elastic relation with unit Young's modulus and Poisson's ratio of 0.3.

According to the range of reference solutions and the curse of dimensionality for three-dimensional data points, only 5 strains equally spacing between $-0.1$ and $0.1$ are sampled for each component, so that the total number of data points is $5^6$. To examine the performance of SLP-UADD algorithm for a noisy data set, the random data sets are generated by adding the Gaussian noise to the precise data points ($\boldsymbol{\varepsilon}_{\text{true}}, \boldsymbol{\sigma}_{\text{true}}$) as: $\boldsymbol{\varepsilon}_{\text{noise}} = \boldsymbol{\varepsilon}_{\text{true}} + \mathcal{N}(0, 0.005)$ and $\boldsymbol{\sigma}_{\text{noise}} = \boldsymbol{\sigma}_{\text{true}} + \mathcal{N}(0, 0.005)$, where $\mathcal{N}(0, 0.005)$ is the Gaussian distribution with the values of mean and variance are equal to 0 and 0.005, respectively.

As listed in Table 10, the numerical parameters for this example are $N_c = 7, L^{(1)} = 1, L_{\min} = 0.2, \rho = 1.5$ and $\text{Tol} = 0.005$, respectively. Since the data set is too sparse, it is set $-0.5 \leq \lambda_i \leq 1.5$ $(i = 1, \dots, N_c)$ to improve the existence of feasible solutions. When there is no feasible solution in a specific iteration, the upper limits of $\lambda_i$ are all added by 1 and the lower limits are all decreased by 1. Once a feasible solution is found, the bounds of $\lambda_i$ are gradually recovered to its initial settings.

The iterative process of $U_{\text{RE}}$ and $\sigma_{\text{RMS}}$ obtained by objective functions of $I(\boldsymbol{U}) =$

$-U_{2271}$ (the $z$-directional displacement of node 757 illustrated by Fig. 14), $I(\boldsymbol{U}) = \boldsymbol{p}^\top \boldsymbol{U}$ and $I(\boldsymbol{U}) = U_{2271}$ are shown in Figs. 15a-c, respectively. It can be found that, the SLP-UADDCM algorithm terminates in 10 iterations for all those cases and the values of $U_{2271}$ converge gradually, which fully demonstrates the effectiveness and robustness of the proposed algorithm under noise. The corresponding values of $U_{2271}$, $U_{\text{RE}}$, $\sigma_{\text{RMS}}$ and the time cost of solving a sequence of linear programming problems are shown in Table 8. Even for such a sparse and noisy data set, the proposed SLP-UADDCM algorithm can still obtain a pair of relatively tight bounds ($\overline{U} = -1.8392$ and $\underline{U} = -1.8698$) of $U_{2271}$ covering its reference value ($U_{\text{ref}} = -1.8650$) and the comparable solution ($U_{\text{s}} = -1.8392$, obtained by setting $I(\boldsymbol{U}) = \boldsymbol{p}^\top \boldsymbol{U}$)). Besides, for all those objective functions, the converged values of $U_{\text{RE}}$ and $\sigma_{\text{RMS}}$ are all smaller than 3%. This clearly illustrates the advantage of the proposed SLP-UADDCM algorithm for dealing with sparse and noisy data sets.

To further investigate the performance of the SLP-UADDCM algorithm, the random data set is refined with a total number of $11^6$ and the variance of Gaussian distribution is set as 0.002. Accordingly, the parameters in Table 10 are set as $L_{\min} = 0.08$ with all others the same. Surprisingly, the algorithm can still converge in 10 iterations even the number of data points is increased by more than 100 times, which demonstrates the effectiveness of the local data point updating strategy. The corresponding solution results such as the upper and lower bounds as well as the solution time costs are also presented in Table 8. Nevertheless, the total solution time is increased by about 100 times in the latter case. This is because, more than 98% of computation time is spent on updating the local data points by using brute-force search algorithm (i.e., searching for the data points closest to vertices of the regular simplex). In order to alleviate this problem, we replace the brute-force search algorithm with a fast approximate nearest neighbor search algorithm implemented in the FLANN-library [33] [34]. The randomized kd-tree algorithm (with 20 random trees) in the FLANN-library is used to update the local data points, and the results are shown in Table 8. It can be found that the solution time for the case where the data set with 1.77 million data points is searched by the approximate nearest neighbor algorithms is actually of the same order of magnitude as that for the case where a data set with much smaller data points is explored by a brute force searching algorithm. In addition, it is also worth noting that in both cases the majority part of the computation time is spent on solving a sequence of linear programming problems during the iteration process. More importantly, the relative errors of displacement and stress obtained by the approximate nearest neighbor search algorithm does not increase, and upper and lower bounds with acceptable accuracy can still be obtained.

Another issue should be mentioned is that, we also encounter some cases, the upper and lower bounds obtained cannot cover the reference solution, as shown in Table 9. This can be understood from the fact that, in the DDCM framework, as long as the number of data points is finite, the problem formulation is non-convex in nature, and thus many local optima exist. Different structural states may be obtained from different initial guesses, and this is also the underlying reason why the classical DDCM algorithms may be trapped by the outliers. Actually, the local convex combination treatment is a local convexification of the original DDCM formulation, and this could improve the robustness of SLP-UADDCM algorithm against noise and outliers to some extent even though the global optimality still cannot be guaranteed. In order to increase the probability of obtaining the global optimal solution, the value of $N_c$ can be increased to enlarge the local convex hulls. The upper and lower bounds of $U_{2271}$ obtained from 30 random data sets with $N_c = 12$ are shown in Fig. 16. It can be observed that for all cases considered, the upper and lower bounds obtained can always cover the reference solution.

*Remark:* The computation time for solving the LPs in the present SLP-UADDCM framework actually constitutes the major part of the total solution time for the considered three-dimensional continuum structure problem when efficient searching algorithm for updating the local data points is adopted. This is quite understandable since the problem of finding the extreme values of structural responses is very time-consuming since it is NP-hard in nature. Furthermore, LP has a relatively simple mathematical structure and can be solved efficiently as compared to other nonlinear programming problems, which are inevitable in classical DDCM framework when nonlinear effects (e.g., geometrical nonlinearity) are considered.

## 5. Concluding remarks

In the present work, a sequential linear programming algorithm (SLP) for the uncertainty analysis-based data-driven computational mechanics (UA-DDCM) is presented. Compared with the existing DDCM paradigm, the distinctive feature of the present approach is that it can provide the upper and lower bounds of the concerned structural response associated with the prescribed data set. In this sense, the present work actually establishes natural links between the three fields of computational mechanics, data science and uncertainty analysis. Numerical examples also demonstrate the effectiveness of the proposed approach.

The present approach also has its limitations and can be improved along different directions. Firstly, although in principle the UA-DDCM formulation does allow for the possibility of obtaining the theoretical upper and/or lower bounds of the concerned structural response, the present solution procedure, however, cannot guarantee that the obtained bounds are the global optimal ones since the current SLP approach does not have the capability of locating the global optimum. Nevertheless, it is worth noting that even though the obtained bounds are not global optimal ones, they are still valuable for evaluating the uncertainties associated with the data set. Obviously, a too large gap between the upper and lower bounds reminds us the necessity of refining the material characterizing process and reducing the uncertainties pertaining to the data set. Therefore, a natural direction for future work is to develop effective approaches which can enhance the possibility of finding the confidence upper and/or lower bounds of the concerned structural response. Secondly, the present algorithm is actually a *multi-point data-based* linear approximation approach which is quite different from the traditional *single-point Taylor's expansion-based algorithm*. Although the numerical experiment clearly demonstrates the effectiveness and robustness of this approach, the corresponding mathematical analysis is still unavailable and needs further exploration. Once the theoretical foundation of this treatment is consolidated, it is expected that the proposed approach can also find applications in the solutions of other types of problems (e.g., structural analysis considering geometrical and material nonlinearities simultaneously, contact analysis), which are difficult to solve by the traditional model-based displacement driven path-following approaches (especially when some bifurcation points exist on the equilibrium path). Some promising results have already been obtained on this aspect and will be reported elsewhere. Finally, for large-scale three-dimensional continuum structures, the number of equality constraints and design variables of SLP-UADDCM algorithm will increase dramatically. In this case, the solution time of the linear programming problems will become the critical factor limiting the efficiency of the proposed algorithm. We also notice that the corresponding optimal solutions of two successive LPs are usually very close especially when the algorithm tends to converge. Under this circumstance, the solution obtained in the previous iteration can be used as the initial solution of the subsequent iteration to accelerate the solution process.

# Appendix: The SLP-UADDCM algorithm for 3D elastic continuum

As illustrated in Table 10, the SLP-UADDCM algorithm for 3D elastic continuum is similar to its 1D counterpart in Table 1. It should be pointed out that, here, all the data points in $\mathcal{D}$ should be sorted according to the algebraic values of the inner product of $\langle \boldsymbol{\varepsilon}_j^{\mathrm{d}}, (1,1,1,1,1,1)^{\top}\rangle$, $j = 1, \dots, N_{\mathrm{d}}$. Moreover, in 3D cases, the third step of updating the data points for local convex hull construction is different from Table 1, which will be explained in detail as follows. Specifically, in order to enhance the feasibility of $\mathcal{P}^{(k)}$ in Eq. (4), the local convex hull in the next step should be able to contain the current stress and strain state as much as possible. However, only selecting a number of data points closest to the current state in the data set will easily cause that the current state is not covered by the local convex hull [3]. An effective treatment for this issue is to cover all possible directions in phase space when local convex hulls are constructed. To this end, we first determine a regular simplex [35] whose vertices are calculated as:

$$\begin{cases} \boldsymbol{\varepsilon}_{1,e}^{\mathrm{s}} = \boldsymbol{\varepsilon}_e^{(k)} - (\gamma_{k+1},\gamma_{k+1},\gamma_{k+1},\gamma_{k+1},\gamma_{k+1},\gamma_{k+1})^{\top}, \\ \boldsymbol{\varepsilon}_{2,e}^{\mathrm{s}} = \boldsymbol{\varepsilon}_{1,e}^{\mathrm{s}} + (p_{k+1},q_{k+1},q_{k+1},q_{k+1},q_{k+1},q_{k+1})^{\top}, \\ \boldsymbol{\varepsilon}_{3,e}^{\mathrm{s}} = \boldsymbol{\varepsilon}_{1,e}^{\mathrm{s}} + (q_{k+1},p_{k+1},q_{k+1},q_{k+1},q_{k+1},q_{k+1})^{\top}, \\ \quad \dots \\ \boldsymbol{\varepsilon}_{7,e}^{\mathrm{s}} = \boldsymbol{\varepsilon}_{1,e}^{\mathrm{s}} + (q_{k+1},q_{k+1},q_{k+1},q_{k+1},q_{k+1},p_{k+1})^{\top}, \end{cases} \quad (8)$$

$$\begin{cases} \boldsymbol{\sigma}_{1,e}^{\mathrm{s}} = \boldsymbol{\sigma}_e^{(k)} - \boldsymbol{\mathcal{D}}(\gamma_{k+1},\gamma_{k+1},\gamma_{k+1},\gamma_{k+1},\gamma_{k+1},\gamma_{k+1})^{\top}, \\ \boldsymbol{\sigma}_{2,e}^{\mathrm{s}} = \boldsymbol{\sigma}_{1,e}^{\mathrm{s}} + \boldsymbol{\mathcal{D}}(p_{k+1},q_{k+1},q_{k+1},q_{k+1},q_{k+1},q_{k+1})^{\top}, \\ \boldsymbol{\sigma}_{3,e}^{\mathrm{s}} = \boldsymbol{\sigma}_{1,e}^{\mathrm{s}} + \boldsymbol{\mathcal{D}}(q_{k+1},p_{k+1},q_{k+1},q_{k+1},q_{k+1},q_{k+1})^{\top}, \\ \quad \dots \\ \boldsymbol{\sigma}_{7,e}^{\mathrm{s}} = \boldsymbol{\sigma}_{1,e}^{\mathrm{s}} + \boldsymbol{\mathcal{D}}(q_{k+1},q_{k+1},q_{k+1},q_{k+1},q_{k+1},p_{k+1})^{\top}, \end{cases} \quad (9)$$

where $\boldsymbol{\mathcal{D}}$ is the scaling matrix[3], and

$$p_{k+1} = \frac{L^{(k+1)}}{6\sqrt{2}}(5+\sqrt{7}), \quad q_{k+1} = \frac{L^{(k+1)}}{6\sqrt{2}}(\sqrt{7}-1), \quad \gamma_{k+1} = \frac{1}{7}(5q_{k+1}+p_{k+1}), \quad (10)$$

respectively. In Eq. (8-10), we choose $N_{\mathrm{c}} = 7$ and the value of $L^{(k+1)}$ represents the Euclidean distance between the regular simplex vertices. Since the vertices $\left((\boldsymbol{\varepsilon}_{j,e}^{\mathrm{s}})^{(k+1)}, (\boldsymbol{\sigma}_{j,e}^{\mathrm{s}})^{(k+1)}\right)$, $j =$ $1, \ldots, N_c$; $e = 1, \ldots, m$ may not coincide with data points, the data pairs $\left((\varepsilon_{j,e}^d)^{(k+1)}, (\sigma_{j,e}^d)^{(k+1)}\right)$, $j = 1, \ldots, N_c$; $e = 1, \ldots, m$ used to construct the local convex hulls in $(k+1)$-th iteration are determined as the data points closest to the vertices of the regular simplex in the data set $\mathcal{D}$ respectively, i.e.,

$$\left\|\left((\varepsilon_{j,e}^d)^{(k+1)}, (\sigma_{j,e}^d)^{(k+1)}\right) - \left((\varepsilon_{j,e}^s)^{(k+1)}, (\sigma_{j,e}^s)^{(k+1)}\right)\right\|$$
$$= \min_{l=1,\ldots,N_d} \left\|(\varepsilon_l^d, \ \sigma_l^d) - \left((\varepsilon_{j,e}^s)^{(k+1)}, (\sigma_{j,e}^s)^{(k+1)}\right)\right\|, \ j = 1, \ldots, N_c; \ e = 1, \ldots, m.$$

(11)

---

[3]The matrix $\boldsymbol{\mathcal{D}}$ guarantees Eqs. (9) and (10) to be consistent since in general the stress and strain components have different magnitudes. In particular, one can use $\boldsymbol{\mathcal{D}} = \mathrm{diag}(D_1, \dots, D_6)$ with $D_i$ denoting the median of the set $\{\sigma_{j,ei}^{\mathrm{d}}/\varepsilon_{j,ei}^{\mathrm{d}}\}$, $j = 1, \dots, N_{\mathrm{c}}$.

This updating strategy of data points for local convex hull construction is illustrated schematically in Fig. 17. Similar to the Algorithm 1 in Table 1, if $\mathcal{P}^{(k)}$ in Eq. (4) is feasible, $L^{(k)}$ is reduced as $L^{(k+1)} = \max(L^{(1)}/\rho^k, L_{\min})$ with $L^{(1)}, \rho, L_{\min}$ denoting the initial length, scaling factor and the low bound of $L^{(k+1)}$. Otherwise, let $L^{(k+1)} = L^{(k)} + 0.1 L^{(k)}$. In addition, the quantity $\lambda_{ej}$, $e = 1, \ldots, m$; $j = 1, \ldots, N_c$ can also be relaxed to increase the feasible region of $\mathcal{P}^{(k+1)}$.

## Acknowledgement

The financial supports from the National Natural Science Foundation (11821202, 11732004, 12002073, 12002077), the National Key Research and Development Plan (2020YFB1709401), 111 Project (B14013) are gratefully acknowledged.

**Figures**

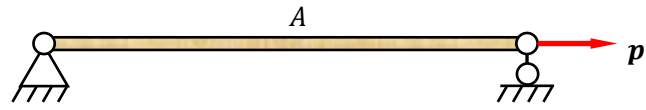

(a)

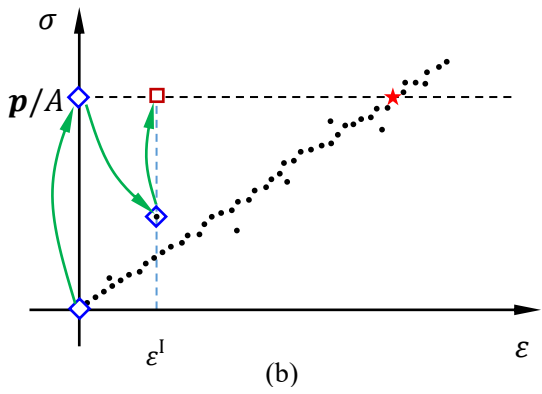

(b)

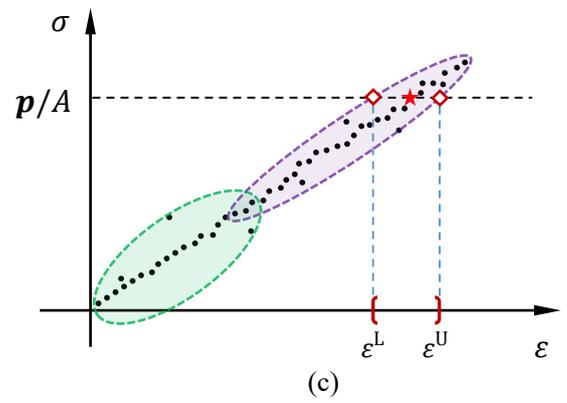

(c)

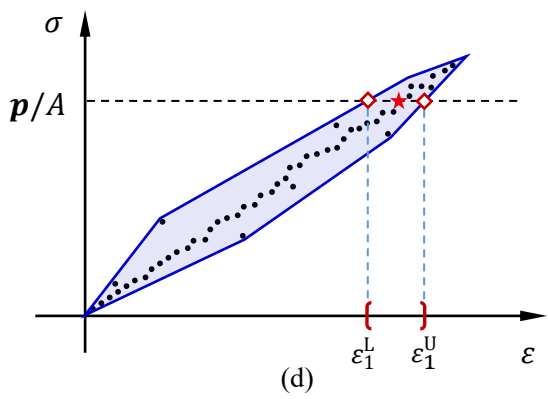

(d)

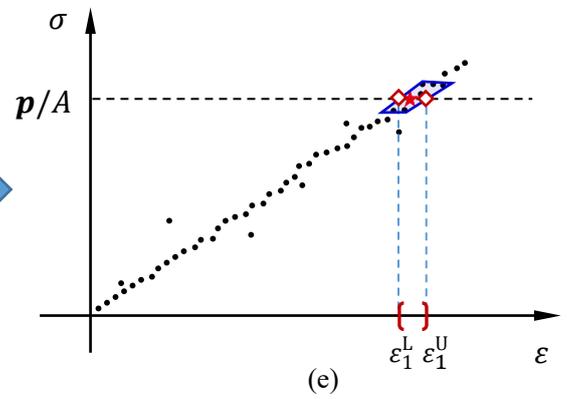

(e)

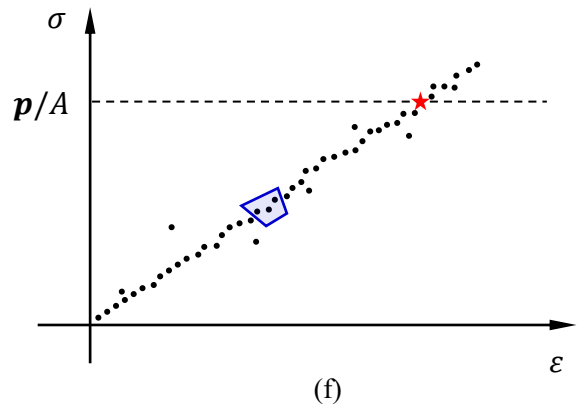

(f)

Fig. 1. (a) An illustrative one-bar example in [5]; (b) an inferior solution $\varepsilon^I$ induced by an outlier in the classical DDCM framework [5]; (c)-(d) confidence bounding intervals $[\varepsilon^L, \varepsilon^U]$ obtained by the UA-DDCM based approach with envelopes constituted by two ellipsoids and a single polygon covering the same data set, respectively; (e) a tighter bounding interval $[\varepsilon^L, \varepsilon^U]$ obtained by the proposed local data set convexification scheme; (f) the case of no feasible solution induced by a possible inappropriate local convexification.

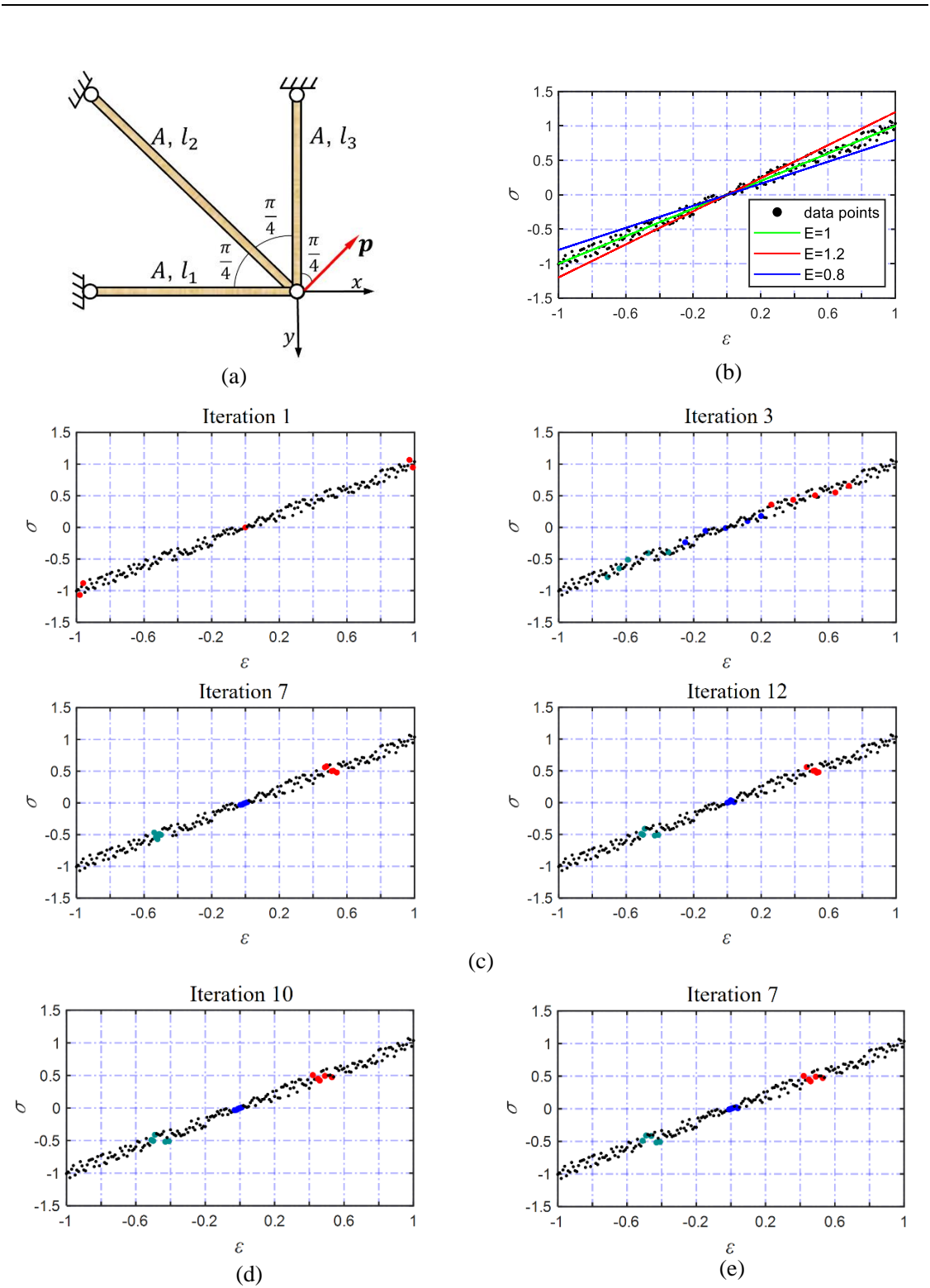

Fig. 2. (a) An illustrative three-bar truss example; (b) the noisy data set; (c) evolution of the local data points associated with each bar involved in convex hull construction at representative iterations with the objective function $I(\boldsymbol{U}) = -U_1$ (red for bar 1, blue for bar 2 and green for bar 3, respectively); (d) converged local data points with the objective function $I(\boldsymbol{U}) = U_1$; (e) converged local data points with the objective function $I(\boldsymbol{U}) = \boldsymbol{p}^\top \boldsymbol{U}$.

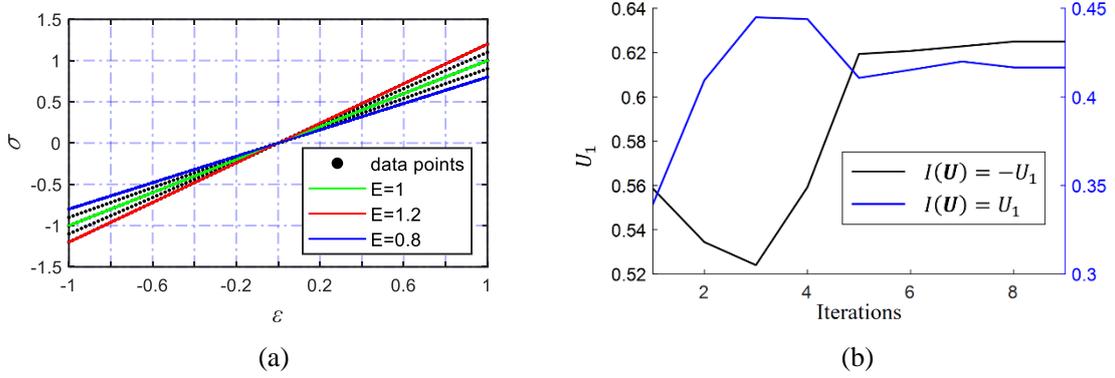

Fig. 3 (a) Regularized data set for three-bar truss example; (b) convergence histories of $U_1$ with different objective functions.

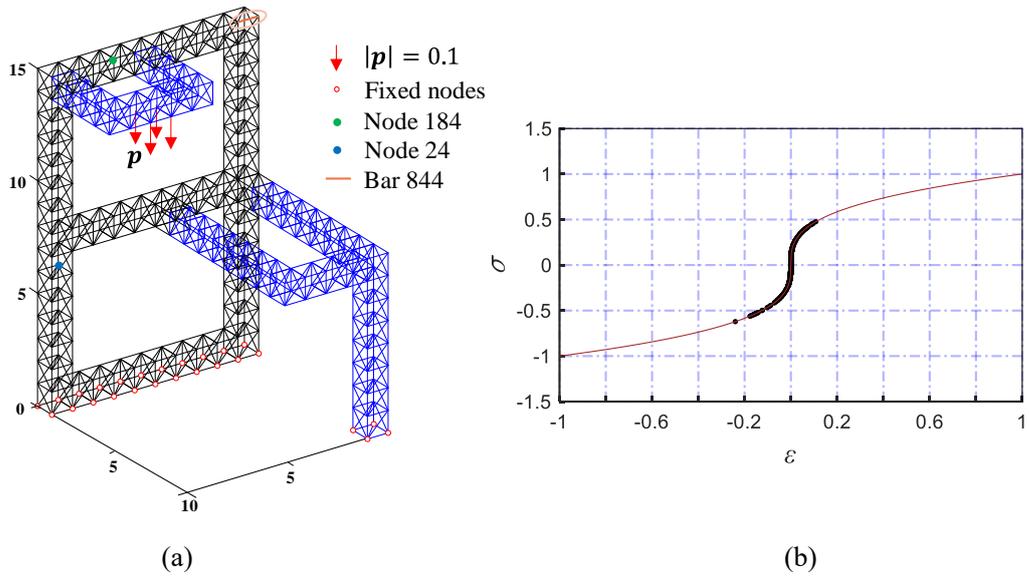

Fig. 4. (a) Problem setting of the three-dimensional truss example; (b) material model ($\sigma = \sqrt[3]{\varepsilon}$) with reference solution values superimposed.

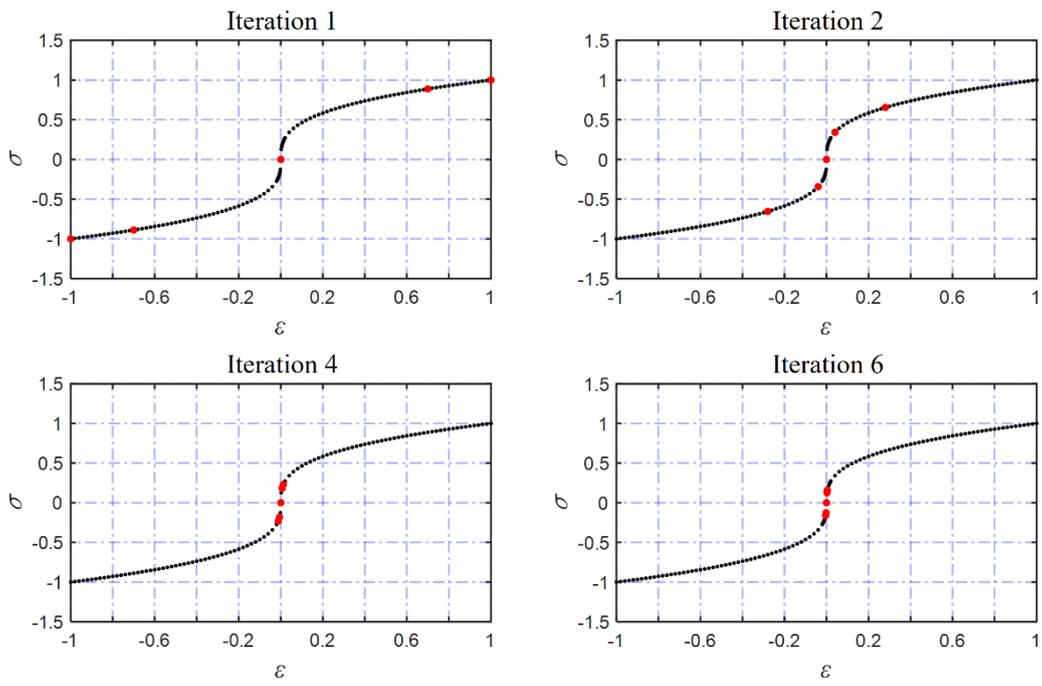

Fig. 5. The local data points associated with the 884-th bar for convex hull construction during the iteration process.

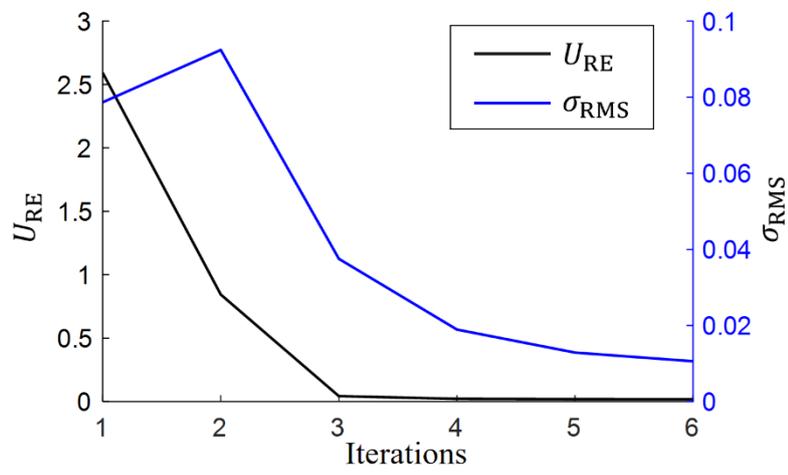

Fig. 6. Convergence histories of $U_{RE}$ and $\sigma_{RMS}$ of the 3D truss example.

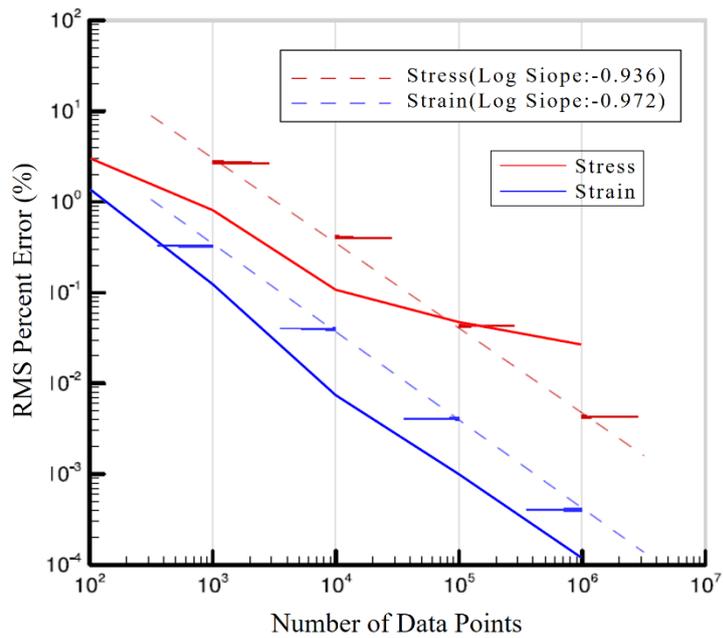

Fig. 7. Convergence of the RMS percent errors of the stress and strain as the increasing of the total number of data points (the solid curves – the present SLP-UADDCM algorithm; the dashed curves – the classical DDCM algorithm in [1]).

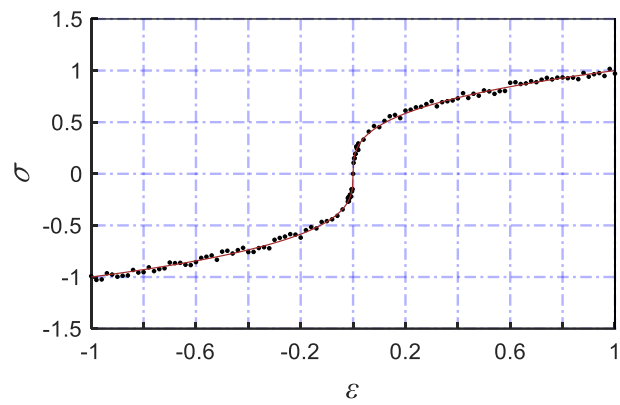

Fig. 8. A random data set containing 121 data points ($\vartheta_0 = 0.04$).

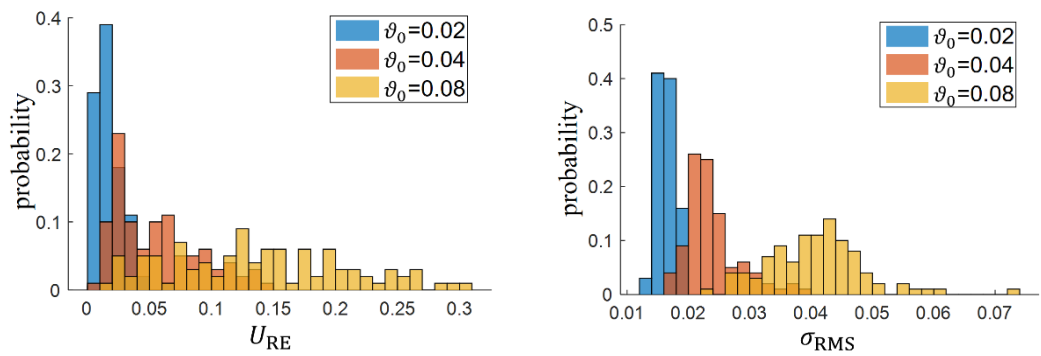

Fig. 9. Histograms of $U_{\text{RE}}$ and $\sigma_{\text{RMS}}$ corresponding to 100 noisy data sets with different values of $\vartheta_0$.

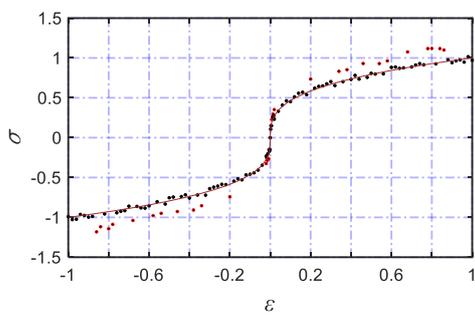 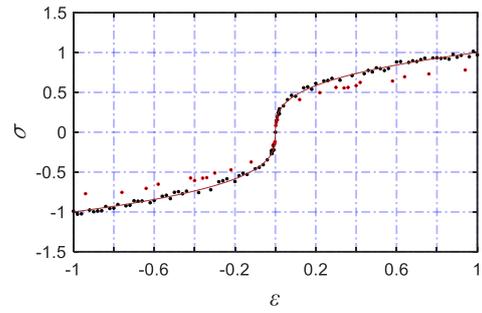

(a) Scaled by 1.2 times of reference values      (b) Scaled by 0.8 times of the reference values

Fig. 10. Two representative data sets with 32 outliers ($N_\mathrm{d} = 121$).

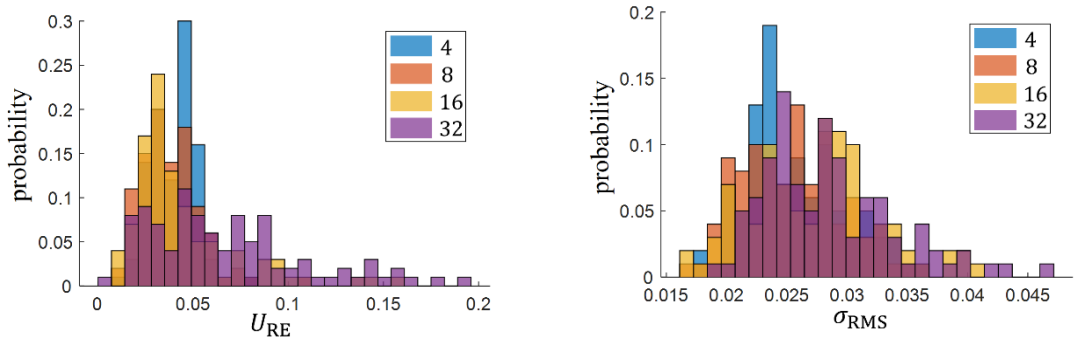

(a) Random outliers scaled by 0.8 times of the reference values

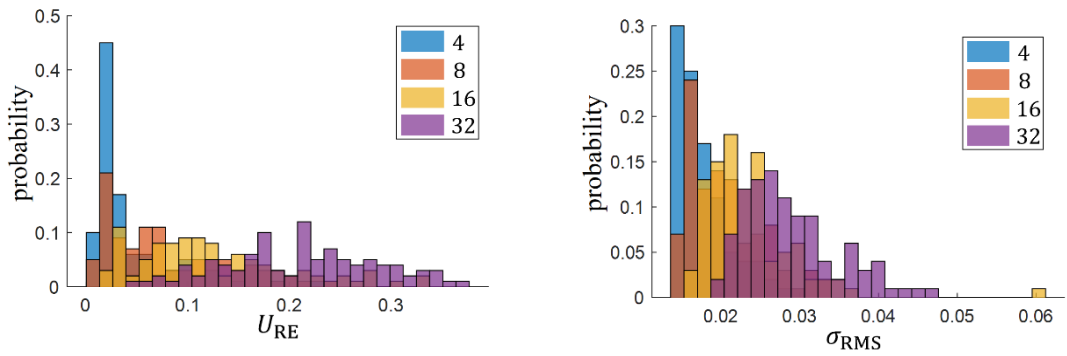

(b) Random outliers scaled by 1.2 times of the reference values

Fig. 11. Histograms of $U_{\text{RE}}$ and $\sigma_{\text{RMS}}$ corresponding to 100 noisy data sets with different numbers of random outliers.

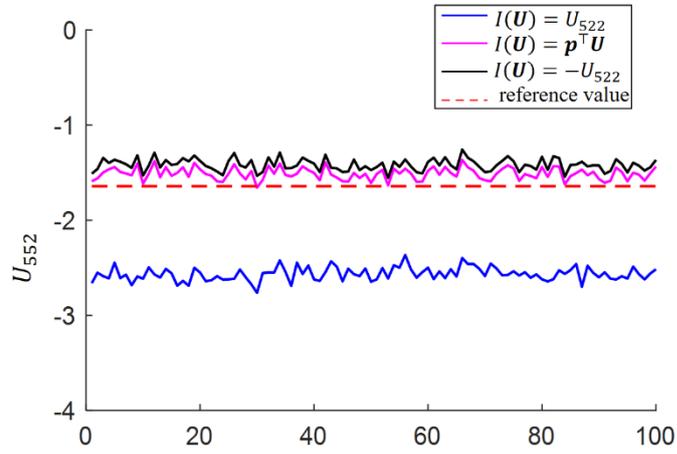

(a) The upper and lower bounds of $U_{552}$

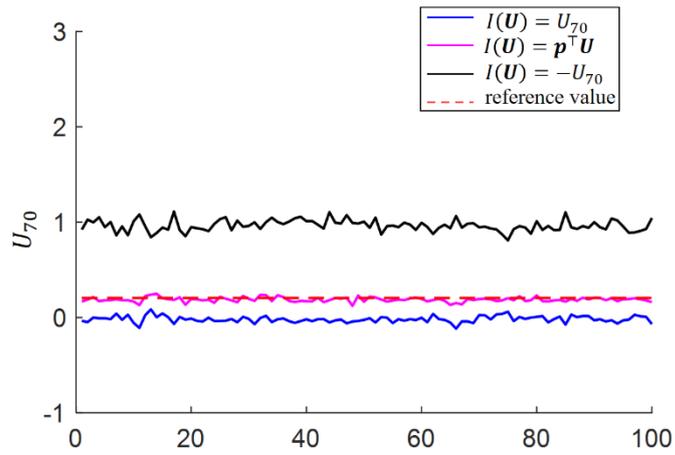

(b) The upper and lower bounds of $U_{70}$

Fig. 12. The upper and lower bounds of $U_{552}$ (a) and $U_{70}$ (b) obtained from 100 random noisy data sets ($\vartheta_0 = 0.04, N_d = 121$).

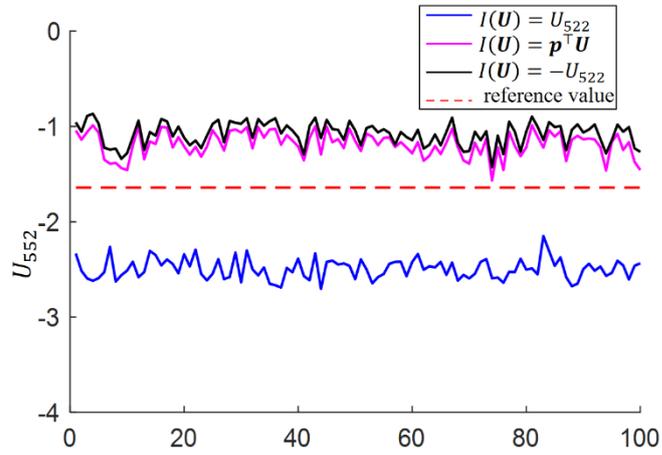

(a) Outliers scaled by 1.2 times of the reference values.

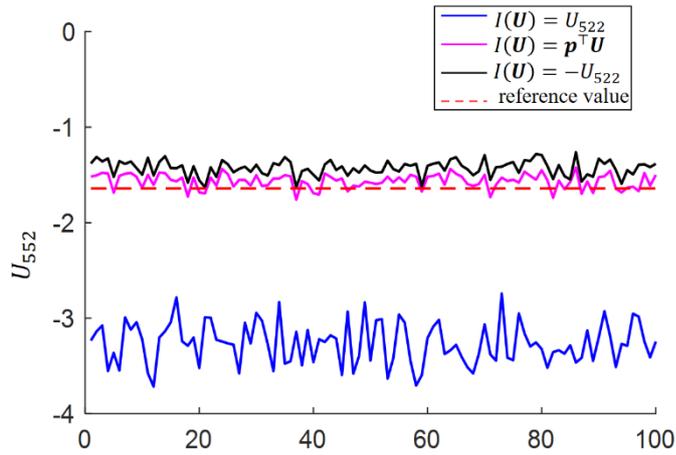

(b) Outliers scaled by 0.8 times of the reference values.

Fig. 13. The upper and lower bounds of $U_{552}$ obtained from 100 random data sets with 16 outliers.

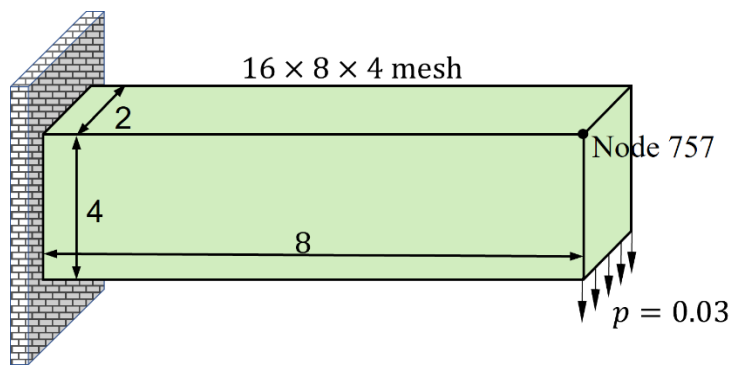

Fig. 14. The problem setting of the 3D cantilever beam example.

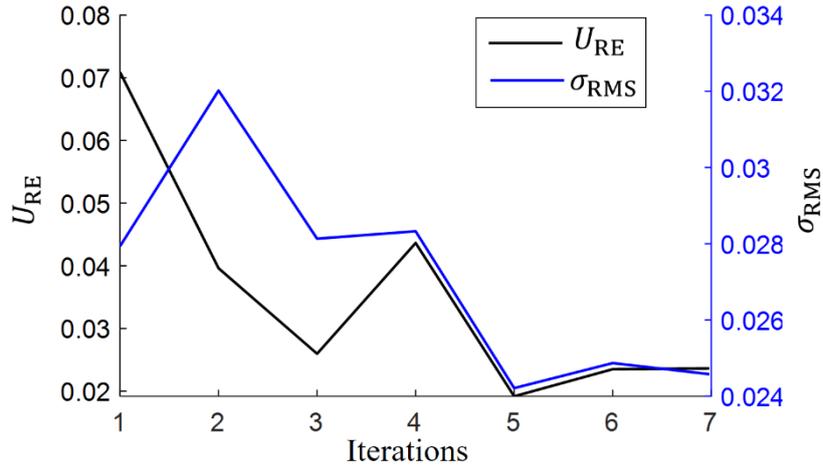

(a) $I(\boldsymbol{U}) = -U_{2271}$

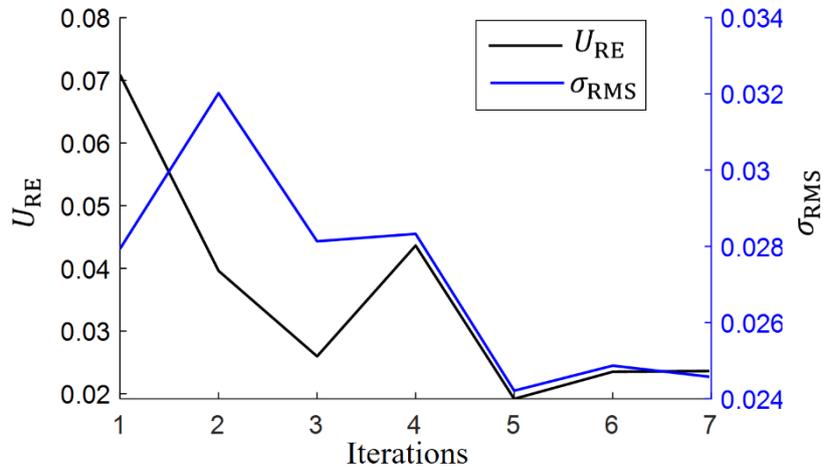

(b) $I(\boldsymbol{U}) = \boldsymbol{p}^\top \boldsymbol{U}$

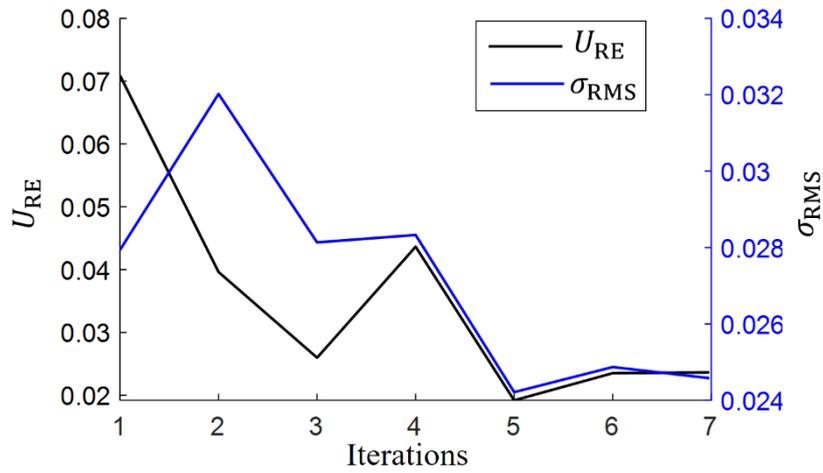

(c) $I(\boldsymbol{U}) = U_{2271}$

Fig. 15. The iteration histories of $U_{\mathrm{RE}}$ and $\sigma_{\mathrm{RMS}}$ for different objective functions ($N_\mathrm{d} = 5^6$).

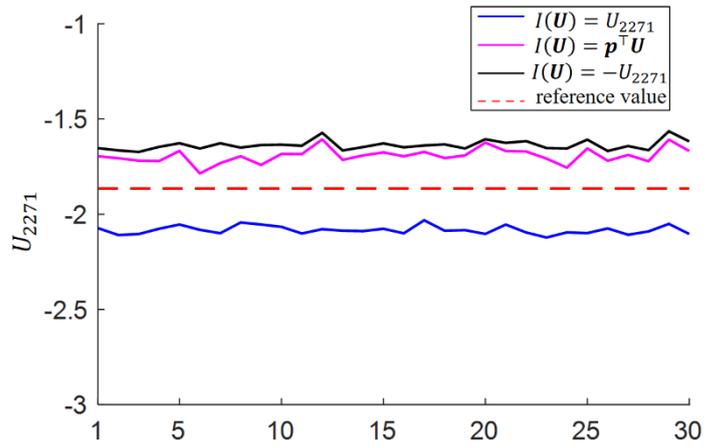

Fig. 16. The upper and lower bounds of $U_{2271}$ obtained from 30 random data sets with $N_c = 12$.

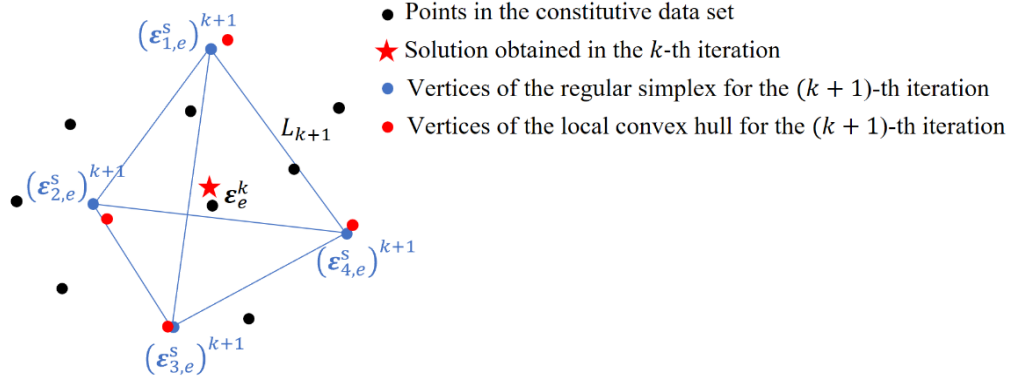

Fig. 17. A schematic illustration of the proposed adaptive local convexification scheme (2D case).

# Tables

Table 1. The SLP algorithm of the UA-DDCM framework for truss structures.

**Input:** Local data sets $\mathcal{D}$, $N_c \geq 3$, $\boldsymbol{b} = (\boldsymbol{b}_1, \dots, \boldsymbol{b}_m)$ and external load $\boldsymbol{p}$.

i) Set $k = 1$ and initialize the data points for local convex hull construction for each bar:

**for all** $e = 1, \dots, m$ **do**

Select $\left\{ \left((\varepsilon_{1,e}^{\mathrm{d}})^1, (\sigma_{1,e}^{\mathrm{d}})^1\right)^{\mathsf{T}}, \dots, \left((\varepsilon_{N_c,e}^{\mathrm{d}})^1, (\sigma_{N_c,e}^{\mathrm{d}})^1\right)^{\mathsf{T}} \right\}$ from $\mathcal{D}$.

**end for**

ii) Find $\boldsymbol{U}^{(k)}$, $\varepsilon_e^{(k)}$, $\sigma_e^{(k)}$:

Solve the linear programming problem $\mathcal{P}^{(k)}$ in Eq. (4).

**If** $\mathcal{P}^{(k)}$ is feasible **then**

$$\begin{pmatrix} \varepsilon_e^{(k)} \\ \sigma_e^{(k)} \end{pmatrix} = \sum_{j=1}^{N_c} \lambda_{ej}^{(k)} \begin{pmatrix} (\varepsilon_{j,e}^{\mathrm{d}})^{(k)} \\ (\sigma_{j,e}^{\mathrm{d}})^{(k)} \end{pmatrix}.$$

**else**

Find $\boldsymbol{U}^{(k)}$, $\varepsilon_e^{(k)}$, $\sigma_e^{(k)}$ following the direct search procedure in classical DDCM algorithm [1].

**end if**

iii) Update the data points used for local convex hull construction for each bar:

**for all** $e = 1, \dots, m$ **do**

Determine $\left\{ \left((\varepsilon_{1,e}^{\mathrm{d}})^{(k+1)}, (\sigma_{1,e}^{\mathrm{d}})^{(k+1)}\right)^{\mathsf{T}}, \dots, \left((\varepsilon_{N_c,e}^{\mathrm{d}})^{(k+1)}, (\sigma_{N_c,e}^{\mathrm{d}})^{(k+1)}\right)^{\mathsf{T}} \right\}$ from $\mathcal{D}$ based on $\left(\varepsilon_e^{(k)}, \sigma_e^{(k)}\right)^{\mathsf{T}}$ following the procedure described below.

**end for**

iv) Check convergence

**If** $\left\| \boldsymbol{U}^{(k)} - \boldsymbol{U}^{(k-1)} \right\|_2 / \left\| \boldsymbol{U}^{(k)} \right\|_2 \leq \mathrm{Tol}$ **then**

$\boldsymbol{U} = \boldsymbol{U}^{(k)}$,

$(\varepsilon_e, \sigma_e)^{\mathsf{T}} = \left(\varepsilon_e^{(k)}, \sigma_e^{(k)}\right)^{\mathsf{T}}$, $e = 1, \dots, m$.

**exit**

**else**

$k = k + 1$, **goto** ii)

**end if**

Table 2. Iteration histories of the concerned nodal displacements obtained with different objective functions for the three-bar truss example.

| Iteration | $I(\boldsymbol{U}) = -U_1$ | | $I(\boldsymbol{U}) = \boldsymbol{p}^\top \boldsymbol{U}$ | | $I(\boldsymbol{U}) = U_1$ | |
|---|---|---|---|---|---|---|
| | $\overline{U}_1$ | $U_2$ | $U_1$ | $U_2$ | $\underline{U}_1$ | $U_2$ |
| 1 | 0.5856 | -0.5487 | 0.3764 | -0.4567 | 0.3764 | -0.4567 |
| 2 | 0.5609 | -0.5665 | 0.4252 | -0.4239 | 0.4252 | -0.4239 |
| 3 | 0.5606 | -0.5701 | 0.4018 | -0.4231 | 0.3813 | -0.4444 |
| 4 | 0.5262 | -0.5226 | 0.4273 | -0.3966 | 0.3908 | -0.4313 |
| 5 | 0.5207 | -0.5123 | 0.3907 | -0.4055 | 0.4640 | -0.4352 |
| 6 | 0.5256 | -0.5383 | 0.4196 | -0.4080 | 0.4750 | -0.4418 |
| 7 | 0.5320 | -0.5106 | 0.4226 | -0.4113 | 0.4780 | -0.4077 |
| 8 | 0.5353 | -0.4749 | | | 0.4282 | -0.4096 |
| 9 | 0.5410 | -0.4993 | | | 0.4184 | -0.4311 |
| 10 | 0.5410 | -0.4625 | | | 0.4184 | -0.4302 |
| 11 | 0.5334 | -0.4629 | | | | |
| 12 | 0.5334 | -0.4629 | | | | |

Table 3. The values of $U_{\text{RE}}$, $\sigma_{\text{RMS}}$, $\varepsilon_{\text{RMS}}$, converged iteration number and the solution time cost with different total numbers of data points ($N_c = 5$).

| $N_d$ | $U_{\text{RE}}$ | $\sigma_{\text{RMS}}$ | $\varepsilon_{\text{RMS}}$ | Iterations | Time(s) |
|---|---|---|---|---|---|
| 41 | 0.1981 | 0.04922 | 0.04261 | 7 | 0.755464 |
| 101 | 0.1188 | 0.03644 | 0.01732 | 7 | 0.783326 |
| 1001 | 0.01356 | 0.01062 | 0.00181 | 9 | 0.932648 |
| 10001 | 5.94E-04 | 0.00159 | 1.30E-04 | 11 | 1.358314 |
| 100001 | 4.79E-05 | 0.000741 | 1.99E-05 | 15 | 2.485592 |

Table 4. The values of $U_{RE}$ and $\sigma_{RMS}$ with different values of $\vartheta_0$ ($N_c = 5$).

| $\vartheta_0$ | $U_{RE}$ | $\sigma_{RMS}$ |
|---|---|---|
| | Mean (variance) | Mean (variance) |
| 0.02 | 1.79% (0.0106) | 1.63% (0.0014) |
| 0.04 | 5.57% (0.0351) | 2.37% (0.0044) |
| 0.08 | 13.82% (0.0714) | 4.07% (0.0080) |

Table 5. The values of $U_{\text{RE}}$ and $\sigma_{\text{RMS}}$ with different numbers of outliers

(scaled by 0.8 times of the reference values).

| Number of outliers | $U_{\text{RE}}$ | $\sigma_{\text{RMS}}$ |
|---|---|---|
|  | Mean (variance) | Mean (variance) |
| 0 | 5.57% (0.0351) | 2.37% (0.0044) |
| 4 | 3.95% (0.0113) | 2.56% (0.0041) |
| 8 | 3.79% (0.0156) | 2.54% (0.0041) |
| 16 | 4.37% (0.0282) | 2.72% (0.0053) |
| 32 | 6.77% (0.0413) | 2.82% (0.0054) |

Table 6. The values of $U_{\text{RE}}$ and $\sigma_{\text{RMS}}$ with different numbers of outliers

(scaled by 1.2 times of the reference values).

| Number of outliers | $U_{\text{RE}}$ | $\sigma_{\text{RMS}}$ |
|---|---|---|
| | Mean (variance) | Mean (variance) |
| 0 | 5.57% (0.0351) | 2.37% (0.0044) |
| 4 | 3.87% (0.0385) | 1.77% (0.0034) |
| 8 | 7.55% (0.0571) | 2.05% (0.0047) |
| 16 | 12.46% (0.0744) | 2.34% (0.0059) |
| 32 | 21.21% (0.0744) | 2.89% (0.0060) |

Table 7. The upper and lower bounds of $U_{552}$ obtained with different total numbers of data points.

| $N_\text{d}$ | $\overline{U}_{552}$ | $U^\text{s}_{552}$ | $\underline{U}_{552}$ |
|---|---|---|---|
| 61 | -1.4727 | -1.6929 | -2.1459 |
| 121 | -1.5868 | -1.6642 | -1.8461 |
| 239 | -1.6112 | -1.6481 | -1.7795 |
| 475 | -1.6264 | -1.6439 | -1.7160 |
| 2363 | -1.6394 | -1.6418 | -1.6473 |

Table 8. The values of $U_{2271}$, $U_{RE}$, $\sigma_{RMS}$ and the solution time cost of the three objective functions with different total numbers of data points (time cost of solving LPs is presented in the brackets).

| | $N_d = 5^6$ (brute-force search method) | | |
|---|---|---|---|
| | $I(U) = -U_{2271}$ | $I(U) = p^\top U$ | $I(U) = U_{2271}$ |
| $U_{2271}$ | -1.8392 | -1.8392 | -1.8698 |
| $U_{RE}$ | 2.36% | 2.36% | 2.36% |
| $\sigma_{RMS}$ | 2.46% | 2.46% | 2.45% |
| Time(s) | 118.13 (103.00) | 118.67 (103.63) | 120.20 (104.56) |
| | $N_d = 11^6$ (brute-force search method) | | |
| | $I(U) = -U_{2271}$ | $I(U) = p^\top U$ | $I(U) = U_{2271}$ |
| $U_{2271}$ | -1.8620 | -1.8632 | -1.8700 |
| $U_{RE}$ | 0.92% | 0.92% | 0.92% |
| $\sigma_{RMS}$ | 1.15% | 1.15% | 1.15% |
| Time(s) | 8443.14 (117.56) | 8454.45 (121.81) | 8613.11 (121.09) |
| | $N_d = 11^6$ (approximate nearest search method) | | |
| | $I(U) = -U_{2271}$ | $I(U) = p^\top U$ | $I(U) = U_{2271}$ |
| $U_{2271}$ | -1.8557 | -1.8557 | -1.8667 |
| $U_{RE}$ | 0.70% | 0.70% | 0.43% |
| $\sigma_{RMS}$ | 1.14% | 1.14% | 1.10% |
| Time(s) | 135.40 (123.24) | 136.03 (124.14) | 135.41 (124.29) |

Table 9. The case where the reference solution ($U_{2271}^{\text{ref}} = -1.8650$) is outside the obtained bounding interval (time cost of solving LPs is presented in the brackets).

| | $N_{\text{d}} = 5^6$ (brute-force search method) | | |
|---|---|---|---|
| | $I(\boldsymbol{U}) = -U_{2271}$ | $I(\boldsymbol{U}) = \boldsymbol{p}^\top \boldsymbol{U}$ | $I(\boldsymbol{U}) = U_{2271}$ |
| $U_{2271}$ | -1.8928 | -1.8928 | -1.9370 |
| $U_{\text{RE}}$ | 3.03% | 3.03% | 3.03% |
| $\sigma_{\text{RMS}}$ | 2.17% | 2.17% | 2.17% |
| Time(s) | 130.06 (112.99) | 130.77 (113.49) | 133.31 (115.65) |
| | $N_{\text{d}} = 11^6$ (brute-force search method) | | |
| | $I(\boldsymbol{U}) = -U_{2271}$ | $I(\boldsymbol{U}) = \boldsymbol{p}^\top \boldsymbol{U}$ | $I(\boldsymbol{U}) = U_{2271}$ |
| $U_{2271}$ | -1.8506 | -1.8506 | -1.8604 |
| $U_{\text{RE}}$ | 0.47% | 0.47% | 0.47% |
| $\sigma_{\text{RMS}}$ | 1.13% | 1.13% | 1.13% |
| Time(s) | 8163.62 (120.45) | 8366.58 (122.08) | 8246.74 (121.13) |
| | $N_{\text{d}} = 11^6$ (approximate nearest search method) | | |
| | $I(\boldsymbol{U}) = -U_{2271}$ | $I(\boldsymbol{U}) = \boldsymbol{p}^\top \boldsymbol{U}$ | $I(\boldsymbol{U}) = U_{2271}$ |
| $U_{2271}$ | -1.8557 | -1.8557 | -1.8642 |
| $U_{\text{RE}}$ | 0.61% | 0.61% | 0.49% |
| $\sigma_{\text{RMS}}$ | 1.10% | 1.10% | 1.10% |
| Time(s) | 118.77 (108.20) | 119.00 (108.25) | 129.90 (118.76) |

Table 10. The SLP-UADDCM algorithm for 3D elastic continuum.

**Input:** Local data sets $\mathcal{D}$, $N_c \geq 7$, strain matrix $\boldsymbol{b}$ and external load vector $\boldsymbol{p}$.

i) Set $k = 1$, and initialize the data points for local convex hull construction for each integration point:

**for all** $e = 1, \ldots, m$ **do**

Choose $\left\{\left(\left(\boldsymbol{\varepsilon}_{1,e}^{\mathrm{d}}\right)^1, \left(\boldsymbol{\sigma}_{1,e}^{\mathrm{d}}\right)^1\right)^{\top}, \ldots, \left(\left(\boldsymbol{\varepsilon}_{N_c,e}^{\mathrm{d}}\right)^1, \left(\boldsymbol{\sigma}_{N_c,e}^{\mathrm{d}}\right)^1\right)^{\top}\right\}$ from $\mathcal{D}$.

**end for**

ii) Find $\boldsymbol{U}^{(k)}, \boldsymbol{\varepsilon}_e^{(k)}, \boldsymbol{\sigma}_e^{(k)}$:

Solve the linear programming problem $\mathcal{P}^{(k)}$ in Eq. (4).

**If** $\mathcal{P}^{(k)}$ is feasible, **then**

$$\begin{pmatrix} \boldsymbol{\varepsilon}_e^{(k)} \\ \boldsymbol{\sigma}_e^{(k)} \end{pmatrix} = \sum_{j=1}^{N_c} \lambda_{ej}^{(k)} \begin{pmatrix} \left(\boldsymbol{\varepsilon}_{j,e}^{\mathrm{d}}\right)^{(k)} \\ \left(\boldsymbol{\sigma}_{j,e}^{\mathrm{d}}\right)^{(k)} \end{pmatrix}.$$

**else**

Find $\boldsymbol{U}^{(k)}, \boldsymbol{\varepsilon}_e^{(k)}, \boldsymbol{\sigma}_e^{(k)}$ following the direct search procedure in the classical DDCM [1].

**end if**

iii) Update the data points used for local convex hull construction for each integration point:

**for all** $e = 1, \ldots, m$ **do**

a) Compute the vertices of each regular simplex $\left(\left(\boldsymbol{\varepsilon}_{j,e}^{\mathrm{s}}\right)^{(k+1)}, \left(\boldsymbol{\sigma}_{j,e}^{\mathrm{s}}\right)^{(k+1)}\right)^{\top}$, $j = 1, \ldots, N_c$ according to $\left(\boldsymbol{\varepsilon}_e^{(k)}, \boldsymbol{\sigma}_e^{(k)}\right)^{\top}$ following Eq. (9) and Eq. (10).

b) Determine $\left\{\left(\left(\boldsymbol{\varepsilon}_{1,e}^{\mathrm{d}}\right)^{(k+1)}, \left(\boldsymbol{\sigma}_{1,e}^{\mathrm{d}}\right)^{(k+1)}\right)^{\top}, \ldots, \left(\left(\boldsymbol{\varepsilon}_{N_c,e}^{\mathrm{d}}\right)^{(k+1)}, \left(\boldsymbol{\sigma}_{N_c,e}^{\mathrm{d}}\right)^{(k+1)}\right)^{\top}\right\}$ from $\mathcal{D}$ based on $\left(\left(\boldsymbol{\varepsilon}_{j,e}^{\mathrm{s}}\right)^{(k+1)}, \left(\boldsymbol{\sigma}_{j,e}^{\mathrm{s}}\right)^{(k+1)}\right)^{\top}$ following Eq. (11).

**end for**

iv) Check convergence

**If** $\left\|\boldsymbol{U}^{(k)} - \boldsymbol{U}^{(k-1)}\right\|_2 / \left\|\boldsymbol{U}^{(k)}\right\|_2 \leq \mathrm{Tol}$ **then**

$\boldsymbol{U} = \boldsymbol{U}^{(k)}$,

$(\boldsymbol{\varepsilon}_e, \boldsymbol{\sigma}_e)^{\top} = \left(\boldsymbol{\varepsilon}_e^{(k)}, \boldsymbol{\sigma}_e^{(k)}\right)^{\top}$, $e = 1, \ldots, m$.

**exit**

**else**

$k = k + 1$, **goto** ii)

**end if**